\documentstyle[amscd,xy]{amsart}

\xyoption{all}\CompileMatrices

\newcommand{\nc}{\newcommand}
\newcommand{\rnc}{\renewcommand}
\nc{\op}[1]{\mathop{{\bf #1}}\nolimits}

\nc{\refeq}[1]{$(\ref{#1})$}

\nc{\sto}[1]{\stackrel{#1}{\longrightarrow}}
\nc{\lbr}{\{\!\{}
\nc{\rbr}{\}\!\}}

\nc{\Hom}{\op{Hom}}
\nc{\Ext}{\op{Ext}}
\nc{\Ker}{\op{Ker}}
\nc{\Coker}{\op{Coker}}
\rnc{\Im}{\op{Im}}
\nc{\codim}{\op{codim}}
\nc{\rk}{\operatorname{rank}\nolimits}
\nc{\can}{\op{can}}
\nc{\id}{\op{id}}
\nc{\ad}{\op{ad}}
\nc{\supp}{\op{supp}}
\nc{\point}{\mathop{{\rm point}}}

\nc{\KK}{{\Bbb K}}
\nc{\NN}{{\Bbb N}}
\nc{\PP}{{\Bbb P}}
\nc{\QQ}{{\Bbb Q}}
\nc{\ZZ}{{\Bbb Z}}

\nc{\BA}{{\bf A}}
\nc{\BB}{{\bf B}}
\nc{\BC}{{\bf C}}
\nc{\BE}{{\bf E}}
\nc{\BH}{{\bf H}}
\nc{\BK}{{\bf K}}
\nc{\BM}{{\bf M}}
\nc{\BN}{{\bf N}}
\nc{\BT}{{\bf T}}
\nc{\BU}{{\bf U}}
\nc{\BX}{{\bf X}}
\nc{\ba}{{\bf a}}
\nc{\bc}{{\bf c}}
\nc{\bi}{{\bf j}}
\nc{\bj}{{\bf i}}
\nc{\bp}{{\bf p}}
\nc{\bq}{{\bf q}}
\nc{\br}{{\bf r}}
\nc{\bk}{{\bf k}}

\nc{\HBX}{{\widehat{\BX}}}
\nc{\HBA}{{\widehat{\BA}}}
\nc{\TBE}{{\widetilde{\BE}}}
\nc{\TBA}{{\widetilde{\BA}}}

\nc{\BRN}{{{\bf NR}}}

\nc{\CE}{{\cal E}}
\nc{\CF}{{\cal F}}
\nc{\CL}{{\cal L}}
\nc{\CN}{{\cal N}}
\nc{\CO}{{\cal O}}
\nc{\CQ}{{\cal Q}}
\nc{\CT}{{\cal T}}
\nc{\CY}{{\cal Y}}

\nc{\FA}{{\frak A}}
\nc{\FE}{{\frak E}}
\nc{\FK}{{\frak K}}
\nc{\FM}{{\frak M}}
\nc{\FU}{{\frak U}}
\nc{\fg}{{\frak g}}
\nc{\fn}{{\frak n}}

\nc{\fsl}{\mathop{\frak s\frak l}}
\nc{\hfgl}{{\widehat{\mathop{\frak g\frak l}}}}
\nc{\hfsl}{{\widehat{\mathop{\frak s\frak l}}}}

\nc{\al}{\alpha}
\nc{\be}{\beta}
\nc{\ga}{\gamma}
\nc{\Ga}{\Gamma}
\nc{\ka}{\kappa}
\nc{\om}{\omega}
\nc{\TH}{\Theta}
\nc{\tth}{{\tilde\theta}}
\nc{\eps}{\varepsilon}
\nc{\tka}{\tilde\kappa}
\nc{\teps}{\tilde\varepsilon}

\nc{\DO}{{\stackrel{{\circ}}{D}}{}}
\nc{\SO}{{\stackrel{{\circ}}{S}}{}}
\nc{\FEO}{{\stackrel{{\circ}}{\FE}}{}}

\nc{\sB}{{\sf B}}
\nc{\sE}{{\sf E}}
\nc{\sK}{{\sf K}}

\nc{\RT}{{\frak R\frak T}}
\nc{\TA}{{\widetilde A}}
\nc{\TD}{{\widetilde D}}
\nc{\TF}{{\widetilde F}}
\nc{\TP}{{\widetilde P}}
\nc{\TR}{{\widetilde R}}
\nc{\TE}{{\widetilde\CE}}
\nc{\TS}{{\widetilde S}}
\nc{\TW}{{\widetilde W}}
\nc{\TGa}{{\widetilde\Ga}}
\nc{\tx}{{\tilde x}}

\nc{\OK}{{\overline{K}}}

\nc{\dx}{\partial}

\rnc{\smile}{\star}
\rnc{\frown}{/}
\nc{\alp}{{\alpha+i}}
\nc{\bep}{{\beta+i}}
\nc{\hp}{h}
\nc{\letimes}{\times_\le}
\nc{\Heis}{\op{Heis}}
\nc{\Diff}{\op{Diff}}

\newsymbol\boxtimes 1202

\newtheorem*{theo}{Theorem}

\newtheorem{th}{Theorem}
\newtheorem{pr}{Proposition}[subsection]
\newtheorem{lm}[pr]{Lemma}
\newtheorem{cor}[pr]{Corollary}

\theoremstyle{definition}

\newtheorem{df}[pr]{Definition}
\newtheorem{rem}[pr]{Remark}

\title[Parabolic sheaves on surfaces and affine Lie algebra $\protect{\hfgl_n}$]
{Parabolic sheaves on surfaces\\and affine Lie algebra $\hfgl_n$}
\author{Michael Finkelberg}
\address{Independent Moscow University, Bolshoj Vlasjevskij pereulok, dom 11,
Moscow 121002 Russia}
\email{fnklberg@@mccme.ru}
\author{Alexander Kuznetsov}
\address{Independent Moscow University, Bolshoj Vlasjevskij pereulok, dom 11,
Moscow 121002 Russia}
\email{sasha@@kuznetsov.mccme.ru}

\begin{document}

\maketitle

\section{Introduction}

The purpose of this paper is to give an example of geometric construction 
(via Hecke correspondences) of certain
representations of the affine Lie algebra $\hfgl_n$. The
construction is similar to the one of \cite{imrn} for the Lie 
algebra $\fsl_n$.

\subsection{The case of $\fsl_n$}

Recall the setup of \cite{imrn}. 
Let $V$ be an $n$-dimensional complex vector space, and let 
$\alpha=(a_1,\ldots,a_{n-1})$ be an $(n-1)$-tuple of nonnegative integers.
We consider $\al$ as a linear combination $\al=\sum a_ii\in\NN[I]$ of simple coroots
$i\in I$ of the Lie algebra $\fsl_n$ via identification $I=\{1,2,\dots,n-1\}$.
Let $\CQ_\al^L$ be 
the space of degree $\al$ quasiflags, that is the space of 
all flags $0\subset E_1\subset\dots\subset E_{n-1}\subset V\otimes\CO_C$
of coherent subsheaves in the trivial vector bundle $V\otimes\CO_C$
over the curve $C=\PP^1$ such that
$$
\operatorname{rank}(E_p)=p,\qquad \deg(E_p)=-a_p.
$$
The space $\CQ^L_\al$ is a smooth compactification of the space $\CQ_\al$
of degree $\al$ maps from the curve $C$ to the flag variety $\BX$ of the 
Lie group~$SL(V)$ (see ~\cite{la}). 
The subspace $\CQ_\al\subset\CQ_\al^L$ is formed by all 
degree $\al$ flags of vector subbundles 
$0\subset E_1\subset\dots\subset E_{n-1}\subset V\otimes\CO_C$.

For $i\in I$ let $\FE_\al^i\subset\CQ^L_\al\times\CQ^L_{\al+i}$ denote the closed
subspace consisting of pairs of flags $(E_\bullet,E'_\bullet)$ such 
that $E'_\bullet\subset E_\bullet$. The subspace 
$\FE_\al^i\subset\CQ^L_\al\times\CQ^L_{\al+i}$ considered as a 
correspondence defines a pair of operators
$$
e_i:H^\bullet(\CQ^L_{\al},\QQ)\to H^{\bullet+2}(\CQ^L_{\al+i},\QQ)
\quad\text{and}\quad
f_i:H^\bullet(\CQ^L_{\al+i},\QQ)\to H^{\bullet-2}(\CQ^L_{\al},\QQ)
$$
Finally, let $h_i:H^\bullet(\CQ^L_{\al},\QQ)\to H^\bullet(\CQ^L_{\al},\QQ)$ 
be the scalar multiplication by $2+2a_i-a_{i-1}-a_{i+1}$.

\begin{theo}[\cite{imrn}] 
The operators $e_i$, $f_i$, $h_i$ provide the vector space 
$$
\BH:=\bigoplus\limits_{\al\in\NN[I]} H^\bullet(\CQ^L_\al,\QQ)
$$
with a structure of $\fsl_n$-module.
\end{theo}

\subsection{Verma submodules}
\label{verma}
Let $\CE_\bullet\in\CQ_{\al_0}\subset\CQ^L_{\al_0}$ be a 
flag of vector subbundles of degree $\al_0$. 
Let $K_\al(\CE_\bullet)\subset\CQ^L_{\al+\al_0}$
denote the closed subspace formed by all quasiflags $E_\bullet$ such that 
$E_\bullet\subset\CE_\bullet$. It is equidimensional of dimension
$|\alpha|=a_1+\ldots+a_{n-1}$. Consider the vector subspace
$$
\BM(\CE_\bullet):=\bigoplus\limits_{\al\in\NN[I]} H^0(K_\al(\CE_\bullet),\QQ)
\subset\BH
$$
spanned by the fundamental cycles of the irreducible components of 
$K_\al(\CE_\bullet)$. It is a $\fsl_n$-submodule of $\BH$, 
isomorphic to the Verma module with
the lowest weight $\al_0+2\rho$.

Let us describe the $\fsl_n$-module structure on $\BM(\CE_\bullet)$ 
in the intrinsic terms of the spaces 
$K_\al(\CE_\bullet)$. We fix the flag $\CE_\bullet$ and
we will write $K_\al$ instead of $K_\al(\CE_\bullet)$ for brevity.
We have to compute
the matrices of the operators $e_i$ and $f_i$ in the bases of fundamental
cycles of top dimensional irreducible components of $K_\al,K_{\al+i}$.
If one viewes $\alpha$ as an element of the coroot lattice of $\fsl_n$,
then $|\al|$-dimensional 
irreducible components of the space $K_\al$ are in one-to-one
corresponedence with Kostant partitions of $\al$ (see ~\cite{mrl}). 
Given a Kostant 
partition $A\in\FK(\al)$ we denote by $\OK_A\subset K_\al$ the
corresponding irreducible component and by $v_A\in H^0(K_\al,\QQ)$
the fundamental class of $\OK_A$.

Let $A\in\FK(\al)$, $A'\in\FK(\al+i)$ be Kostant partitions.
In order to compute the matrix coefficient $\eps_i(A,A')$ of
the operator $e_i$ we should describe the intersection
$$
\FE_A^{A'}:=\FE_\al^i\cap p^{-1}(\OK_A)\cap q^{-1}(\OK_{A'}),
$$
where $p:\FE_\al^i\to\CQ^L_\al$ and $q:\FE_\al^i\to\CQ^L_{\al+i}$
are the projections. 
More precisely, we need to know all irreducible components of $\FE_A^{A'}$
which are dominant over $\OK_{A'}$. One can check that all these components
have the expected dimension and the 
intersection $\FE_\al^i\cap p^{-1}(\OK_A)$
is transversal along them. Hence the matrix coefficient $\eps_i(A,A')$
is equal to the sum of degrees of these components over $\OK_{A'}$.

Similarly, in order to compute the matrix coefficient $\phi_i(A',A)$
of the operator $f_i$ we need to know all irreducible components 
of $\FE_A^{A'}$ which are dominant over $\OK_{A}$. However, in contrast
to the case of $e_i$, the situation here is rather complicated.
Namely, in some cases these irreducible components have the expected
dimension and then the matrix coefficient is equal to the degree of 
these components over~$\OK_A$. But in some cases the dimension of 
these irreducible components exceeds the expected dimension by 1
(the excess intersection). In these cases we also need to describe
the excess intersection line bundle on these components. Then the 
matrix coefficient is equal to the sum of the degrees of the restrictions 
of this line bundle to a generic fiber of these components over $\OK_A$.

\subsection{The case of $\hfgl_n$: a wishful thinking} 

Let $C\subset S$ be a smooth compact 
curve of genus $g=g(C)$ in a smooth compact surface $S$. Let us fix a sequence 
$\gamma=(c_i),\ i\in\ZZ$, of cohomology classes in $H^2(S,\ZZ)$, 
such that $c_{i+1}=c_i+[C]$, and a sequence $\al=(a_i),\ i\in\ZZ$, 
of cohomology classes in $H^4(S,\ZZ)=\ZZ$ such that 
$a_{i+n}=a_i+c_i\cdot[C]+n[C]^2/2$. Recall that a parabolic sheaf of
rank $n$ on $(S,C)$ is an infinite flag of torsion free rank $n$ coherent
sheaves $\ldots\subset E_{-1}\subset E_0\subset E_1\subset\ldots$ such that
$E_{i+n}=E_i(C)$. Let $\CY(\gamma,\alpha)$ be the moduli space of $\mu$-stable
rank $n$ parabolic sheaves such that $ch_1(E_i)=c_i,\ ch_2(E_i)=a_i$ (see ~\cite{y}).
For an $n$-periodic sequence of nonnegative integers $\beta=(b_i)$, 
let $\FE^\beta(\gamma,\alpha)\subset\CY(\gamma,\alpha)\times
\CY(\gamma,\alpha+\beta)$ be the closed subspace formed by all pairs of parabolic sheaves
$(E_\bullet,E'_\bullet)$ such that $E'_\bullet\subset E_\bullet$, and 
$E_\bullet|_{S-C}=E'_\bullet|_{S-C}$. 

We expect that $\FE^\beta(\gamma,\alpha)$
is equidimensional of dimension $\dim\CY(\gamma,\alpha)+|\beta|$
(where $|\beta|=\sum_{i=0}^{n-1}b_i$), and its irreducible components are
parametrized by the set $\FK(\beta)$ of
isomorphism classes of $\beta$-dimensional nilpotent
representations of the cyclic quiver with $n$ vertices (affine quiver of type
$\TBA_{n-1}$). For an isomorphism class 
$\kappa\in\FK(\beta)$ let us denote the 
corresponding irreducible component of $\FE^\beta(\gamma,\alpha)$ by
$\FE_\kappa^\beta(\gamma,\alpha)$.
The set of the above isomorphism classes over all
$\beta\in\NN[\ZZ/n\ZZ]$ forms a basis of the~{\em Hall algebra} $\BH_n$ of
the category of nilpotent representations of the cyclic quiver~$\TBA_{n-1}$.
An irreducible component $\FE_\kappa^\beta(\gamma,\alpha)$ viewed as a
correspondence between $\CY(\gamma,\alpha)$ and $\CY(\gamma,\alpha+\beta)$
defines the map $H^\bullet(\CY(\gamma,\alpha),\QQ)\to
H^\bullet(\CY(\gamma,\alpha+\beta),\QQ)$, and we expect that this way one obtains
the action of $\BH_n$ on $\oplus_\alpha H^\bullet(\CY(\gamma,\alpha),\QQ)$.

According to ~\cite{s}, $\BH_n$ is isomorphic to the positive part 
of the enveloping
algebra~$U(\hfgl_n)$. We expect that the above correspondences transposed
define the action of the negative part of $U(\hfgl_n)$, and together they
generate the action of $U(\hfgl_n)$. 

\subsection{Back to Verma}

At the moment we are unable to carry out the above program. We have to 
restrict ourselves to the affine analog of ~\ref{verma}. Namely, we fix a
parabolic sheaf consisting of {\em locally free sheaves}
$\ldots\subset\CE_{-1}\subset\CE_0\subset\CE_1\subset\ldots$
For an $n$-periodic sequence of nonnegative integers $\alpha=(a_i)$ we define
$K_\alpha=K_\alpha(\CE_\bullet)$ as the space formed by all parabolic sheaves
$E_\bullet$ such that $E_\bullet\subset\CE_\bullet$, and $\CE_i/E_i$ is
concentrated on $C\subset S$ and has (finite) length $a_i$ for any $i\in\ZZ$.
We prove that $K_\alpha$ is equidimensional of dimension $|\alpha|$, and its
irreducible components are naturally parametrized by the set $\FK(\alpha)$
of isomorphism classes of nilpotent representations of dimension $\alpha$ 
of the cyclic quiver $\TBA_{n-1}$. The fundamental classes of these components form a
basis of $\BM=\bigoplus\limits_\alpha H^0(K_\alpha,\QQ)$. For $\kappa\in\FK(\alpha)$ we
denote the corresponding basis element by $v_\kappa$. 

For $i\in\ZZ/n\ZZ$ let $\alp$ denote the sequence $(a'_j)$ such that
$a'_j=a_j+\delta_{ij}$. Let $\FE^i_\alpha\subset K_\alpha\times K_{\alp}$
be the correspondence formed by all the pairs $(E_\bullet,E'_\bullet)$ such
that $E'_\bullet\subset E_\bullet$.  It defines
the maps 
$$
e_i:\ H^0(K_\alpha,\QQ)\to H^0(K_{\alp},\QQ),\quad\text{and}\quad
f_i:\ H^0(K_{\alp},\QQ)\to H^0(K_\alpha,\QQ),
$$ 
the matrix coefficients in the basis $v_\kappa$ are defined 
similarly to ~\ref{verma}. We prove that
$e_i,f_i,\ i\in\ZZ/n\ZZ$, generate the action of $\hfsl_n$ on $\BM$.

The (restricted) dual space $\BM^*$ may be identified with the polynomial
algebra $\QQ[x_\theta]$ 
on infinitely many generators parametrized by the indecomposable
nilpotent representations of $\TBA_{n-1}$. The dual action of
(Chevalley generators of) $\hfsl_n$ on $\BM^*$ is realized by the explicit
first order differential operators in the coordinates $x_\theta$.

On the other hand, $e_i\in\hfsl_n,\ i\in\ZZ/n\ZZ$, generate the positive
part of the universal enveloping algebra $U^+(\hfsl_n)$ which is naturally
embedded into $\BH_n$ (as the subalgebra generated by the isomorphism classes
of indecomposable nilpotent representations). We write down explicit formulae
for the action of $U(\hfgl_n)\supset U^+(\hfgl_n)=\BH_n\supset U^+(\hfsl_n)$
by differential operators in the coordinates $x_\theta$. At the moment we
cannot prove the geometric meaning behind these formulae (see though the
Remark ~\ref{meaning}). Let us only mention
that the central charge of $\BM$ equals $(2-2g(C))n+[C]^2$ (recall that 
$[C]^2$ equals the degree of the normal bundle $\CN_{C/S}$).

\subsection{} 

Let us say a few words about the structure of the paper.
In \S2 we recall the necessary information about the nilpotent 
representations of cyclic quivers in various categories. In \S3 we
study the space $K_\alpha$ together with its projection to a configuration
space of $C$. We prove that all the fibers of this projection admit a cell
decomposition, and compute dimensions of all the cells. In \S4 we
study the correspondence $\FE^i_\alpha\subset K_\alpha\times K_{\alp}$,
and describe its irreducible components dominant over the topdimensional
components of $K_\alpha$ and $K_{\alp}$. It appears that for every
component of $\FE^i_\alpha$ its projection to $K_{\alp}$ is
semismall. In \S5 we define geometrically the matrix coefficients of
the operators $e_i,f_i,h_i\in\hfsl_n$ in the basis of topdimensional
components of $K_\alpha,\ \alpha\in\NN[\ZZ/n\ZZ]$, and compute them explicitly.
In \S6 we realize the $\hfsl_n$-module $\BN$ dual to $\BM$ of \S5
in the polynomial functions on an infinite-dimensional affine space. The action
of $\hfsl_n$ on $\BN$ is given by explicit differential operators. This
realization is similar to Kostant's construction of a dual Verma module over
a semisimple Lie algebra in the sections of a line bundle over the big Schubert
cell. Finally, in \S7 we write down explicit differential operators
extending the $\hfsl_n$-action on $\BN$ to $\hfgl_n$. It is likely that the
resulting $\hfgl_n$-module is a contragredient Verma module.

\subsection{} 

The work ~\cite{imrn} might be viewed as a
globalization of the geometric construction of $U(\fsl_n)$ discovered 
in ~\cite{blm}, ~\cite{g} (from a nilpotent neighbourhood of a point $x\in C$
to the global curve $C$). Similarly, the present work may be viewed as a 
globalization of the geometric construction of $U(\hfsl_n)$ discovered 
in ~\cite{gv}, ~\cite{lu}. In another direction, ~\cite{imrn} was generalized 
in ~\cite{ffkm} to arbitrary simple finite dimensional algebras. It would be
extremely important to generalize the present work along these lines. 

Quite naturally, in this paper $n\geq2$. But in fact, many considerations
below make sense for $n=1$ when the moduli space of parabolic sheaves becomes
just the punctual Hilbert scheme of $S$. We were strongly influenced by the
beautiful book ~\cite{n}, especially chapters 7--9 (see also ~\cite{gr}).
In another direction,
we were strongly motivated by the suggestion of V.~Ginzburg back in 1997 to
study the Drinfeld compactification of the space of maps from $\PP^1$ to an
affine Grassmannian. We learned from him about the parallel between Laumon
and Drinfeld compactifications on the one hand, and Gieseker and Uhlenbeck
compactifications, on the other.  

We would like to thank V.~Baranovsky for bringing the 
reference ~\cite{y} to our attention, and E.~Vasserot for sending us the
preprint ~\cite{s}. We are grateful to R.~Bezrukavnikov for an important 
encouragement and advice,
and to B.~Bakalov and T.~Pantev for useful discussions. 

While this paper was written, the first author enjoyed the hospitality of 
the Insitute for Advanced Study and the support of the NSF grant DMS 97-29992,
and the second author enjoyed the hospitality of the Institut des Hautes
\'Etudes Scientifiques.

\section{Notations}





\subsection{Nilpotent representations of the cyclic quiver}

Let $\TBA_{n-1}$ denote the cyclic quiver with $n$ vertices and $\BA_\infty$
denote the infinite linear quiver. A representation $M$ of the quiver 
$\BA_\infty$ is a $\ZZ$-graded vector space with an operator $A$ of degree~$1$.
A~representation of the quiver $\TBA_{n-1}$ is the same as a $n$-periodic
representation of the quiver $\BA_\infty$, that is a representation $M$
with an isomorphism $M[n]\cong M$, where $[n]$ is the functor shifting 
the grading by $-n$.

A representation $M$ is called {\em nilpotent} if $A^N=0$ for $N\gg0$.
Let $\BRN$ denote the category of nilpotent representations of the 
quiver $\BA_\infty$ and $\BRN_n$ the category of nilpotent representations
of the quiver $\TBA_{n-1}$. 

The dimension of a representation $M$ of the quiver $\BA_\infty$
is just a sequence of nonnegative integers equal to the dimensions of
the graded components of $M$. If~$M$ is a representation of
the quiver $\TBA_{n-1}$ then $\dim M$ is a $n$-periodic sequence. 
 
Recall the well known classification of indecomposable objects of $\BRN_n$.
They are classified up to isomorphism by pairs $(p,q)$ where $p\leq q$ are
integers defined up to simultaneous translation by a multiple of $n$. The
representation corresponding to $(p,q)$ is denoted by $M_{(p,q)}$. It has a
basis $e_p,e_{p+1},\ldots,e_q$ with $e_j$ of degree $j\pmod{n}$, and we have
$e_p\to e_{p+1}\to\ldots\to e_q\to0$ in the representation. We will denote
the Grothendieck group $\BK(\BRN_n)$ by $\BK_n$.
It has a basis $[M_{(p,q)}],\ p,q$ as above. We will denote this basis
by $R^+$. An element of $R^+$ will be called {\em raiz}, and sometimes
a raiz $[M_{(p,q)}]$ will be denoted simply by $(p,q)$. 
Given an integer $s$ we can identify the set $R^+$ with the set
$$
R^+_s=\{(p,q)\ |\ \text{$p\le q$ and $s\le q\le s+n-1$}\}
$$
A raiz $(i,i)$ will be called {\em simple}, and sometimes will be denoted
simply by $i$. We denote by $I\cong\ZZ/n\ZZ$ the set of all simple raiz.

\begin{df}
We say that a raiz $\theta=(p,q)\in R^+$ begins (resp. ends) at a simple 
raiz $i$ iff $i=p \bmod n$ (resp. $i=q \bmod n$). We denote by
$\sB_i\subset R^+\cup\{0\}$ 
(resp. $\sE_i\subset R^+\cup\{0\}$)
the set consisting of $0$ and of all raiz beginning (resp. ending) at $i$.
\end{df}

\begin{df} Given two indecomposable representations $\theta,\vartheta$
such that for some $i$ we have $\theta\in\sB_i,\ 
\vartheta\in\sE_{i-1}$, there is a unique (isomorphism class of an)
indecomposable representation $\eta$ fitting into exact sequence
$$0\to\theta\to\eta\to\vartheta\to0$$ We say that $\eta=\vartheta\smile\theta,\
\vartheta=\eta\frown\theta$.
\end{df}

The dimension of a representation may be viewed as an element
of the lattice~$Y$ of $n$-periodic sequences of integers.
We identify a simple raiz $i\in I$ with the dimension of the
corresponding simple representation, thus if $\al\in Y$
then $\alp\in Y$ is the same as $\al+\dim i$. 

The above identification gives rise to the identification $Y=\ZZ[I]$. 
Thus the dimension may be viewed as a map
$$
\dim:\BK_n\to \NN[I]\subset\ZZ[I]=Y.
$$


We will consider the following elements of the dual lattice $Y^\vee$.
For $i\in\ZZ/n\ZZ$ we define 
\begin{equation}\label{iprime}
\langle i',y_\bullet\rangle=2y_i-y_{i-1}-y_{i+1}
\end{equation}
where $y_\bullet=(\ldots,y_{-1},y_0,y_1,\ldots)\in Y$. 

The $I\times I$ matrix
$$
a_{ij}=\langle i', j\rangle=\begin{cases}
 2, & \text{if $j=i \bmod n$;}\\
-1, & \text{if $j=i\pm1 \bmod n$ and $n\ne 2$;}\\
-2, & \text{if $j=i+1 \bmod n$ and $n=2$;}\\
 0, & \text{if $j\ne i,\ i\pm1 \bmod n$.}
 \end{cases}
$$
is the affine Cartan matrix of type $\TBA_{n-1}$.

For $\alpha\in\NN[I]$ we define $|\alpha|$ as the sum of all
coordinates of $\alpha$. For $\alpha,\beta\in\NN[I]$ we say that 
$\alpha\leq\beta$ iff $\beta-\alpha\in\NN[I]$.

Consider also the lattice $Y^{(2)}\supset Y$ of sequences of 
integers $n$-periodic modulo linear term, that is
$$
Y^{(2)}=\{y_\bullet\ |\ \text{$y_{p+n}-y_p=ap+b$ for some $a,b\in\ZZ$}\}.
$$
Also, we define
$$
\rho=(\dots,\rho_{-1},\rho_0,\rho_1,\dots)\in Y^{(2)}\otimes\QQ,\qquad \rho_p=-p^2/2.
$$

\begin{lm}
The functional $i'$ defined by \refeq{iprime} is well defined
on the lattice $Y^{(2)}$ $($that is $\langle i',y_\bullet\rangle$ depends
only on $i \bmod n$ when $y_\bullet\in Y^{(2)})$. Moreover, we have
$$
\langle i',\rho\rangle=1
$$
for all $i\in I$.
\end{lm}
\begin{pf}
Evident.
\end{pf}




\subsection{Partitions}

Assume that we have a set $X$ and a function $\xi:X\to(\NN[I]-\{0\})$.
For any $\al\in\NN[I]$ we define an ordered $m$-terms 
partition $\TA$ of $\al$ with respect to $(X,\xi)$ as 
a map $\TA:[m]=\{1,\dots,m\}\to X$ such that 
$$
\sum_{p=1}^m\xi(\TA(p))=\al.
$$
We denote the set of all ordered $m$-terms partitions of $\al$ 
with respect to $(X,\xi)$ by $\TP_{X,\xi}^m(\al)$. The group
of permutations $S_m$ acts naturally on the set $\TP_{X,\xi}^m(\al)$. We denote
by $S_\TA\subset S_m$ the stabilizer subgroup of $\TA\in\TP_{X,\xi}^m(\al)$.

We define an unordered $m$-terms partition (or, more simply, a partition) 
$A$ of $\al$ with respect to $(X,\xi)$ as an $S_m$-orbit in the set 
$\TP_{X,\xi}^m(\al)$, and denote by 
$$
P_{X,\xi}(\al)=\bigsqcup_{m=0}^\infty P_{X,\xi}^m(\al)=
\bigsqcup_{m=0}^\infty\TP_{X,\xi}^m(\al)/S_m
$$ 
the set of all partitions of $\al$ with respect to $(X,\xi)$.

Given a partition $A$ we denote by $\TA$ its {\em ordering},
that is any representative of $A$ in the set of ordered partitions.
Let $\al_p=\TA(p)$ ($p=1,\dots,m$). We will denote the partition 
$A$ by $\lbr\al_1,\dots,\al_m\rbr$, and for $A=\lbr\al_1,\dots,\al_m\rbr$ 
we will denote
$$
|A|=\xi(\al_1)+\dots+\xi(\al_m)\in\NN[I],\qquad \sK(A)=m,
$$
and
$$
m(\al,A)=\#\{p\in[m]\ |\ \al_p=\al\}.
$$
Thus, we always have $A\in P_{X,\xi}(|A|)$.

We will use the following types of partitions.

\begin{description}

\item[Usual partitions] Here we put $X=\NN[I]-\{0\}$, $\xi=\id$. 
We denote the set of usual partitions of $\al$ by $\Ga(\al)$.

\item[Kostant partitions] Here we put $X=R^+$, $\xi=\dim$. 
We denote the set of Kostant partitions of $\al$ by $\FK(\al)$.

\item[Multipartitions] Here we put 
$X=\bigsqcup_{\ga\in(\NN[I]-\{0\})}\FK(\ga)$,
$\xi(\ka)=|\kappa|$. 
We denote the set of multipartitions of $\al$ by $\FM(\al)$.

\end{description}

Note that if $\theta\in R^+$ is a raiz,
then the set $\FK(\dim\theta)$ of Kostant partitions of $\dim\theta$ contains 
an element $\lbr\theta\rbr$. Such Kostant partition is called
{\em a simple Kostant partition}. 
A multipartition $\mu=\lbr\ka_1,\dots,\ka_m\rbr$
is called {\em a simple multipartition} if all $\ka_p$ are simple Kostant
partitions.

We have the following natural maps: $\dim:\FK(\al)\to\Ga(\al)$,
$\FK(\al)\hookrightarrow\FM(\al)$ (the set of Kostant partitions of $\al$
is identified with the set of simple multipartitions of $\al$),
and the projection $|\ |:\FM(\al)\to\Ga(\al)$
(we define $|\lbr\ka_1,\dots,\ka_m\rbr|:=\lbr|\ka_1|,\dots,|\ka_m|\rbr$).

\subsection{Configuration spaces}

Let $C^\al$ denote the configuration space of effective divisors
on the curve $C$ with coefficients in $\NN[I]$ of degree $\al$.
If $\al=\sum_{i\in I}a_ii$ then the space $C^\al$ is isomorphic
to the product of the symmetric powers of $C$, more presicely
$$
C^\al\cong \prod_{i\in I}S^{a_i}C.
$$

The space $C^\al$ carries the natural diagonal stratification, the strata
of which are in one-to-one correspondence with partitions of $\al$:
$$
C^\al=\bigsqcup_{\Ga\in\Ga(\al)}C^\al_\Ga.
$$
The stratum $C^\al_\Ga$, corresponding to the ordered partition 
$\TGa=\lbr\al_1,\dots,\al_m\rbr$ consists of all divisors 
$\sum_{p=1}^m\al_px_p$, where $x_1,\dots,x_p$ are pairwise 
distinct points of $C$. Thus we have an isomorphism
$$
C^\al_\Ga\cong(C^m-\Delta)/S_\TGa,
$$
where $\Delta\subset C^m$ is the big diagonal.

If $\mu\in\FM(\al)$ is a multipartition and $\Ga=|\mu|\in\Ga(\al)$
then it is clear that $S_{\tilde\mu}\subset S_\TGa$. Let $C^\al_\mu$ 
denote the quotient $(C^m-\Delta)/S_{\tilde\mu}$. The space $C^\al_\mu$ 
is a ${\displaystyle\frac{|S_\TGa|}{|S_{\tilde\mu}|}}$-fold
covering of the space $C^\al_\Ga$.

\subsection{Nilpotent $\TBA_{n-1}$-modules over $C$}
Let $\al\in\NN[I]$.
Recall that the isomorphism classes of $\al$-dimensional
objects of the category $\BRN_n$ are in one-to-one correspondence with
Kostant partitions of $\al$. We denote by $\ka(M)\in \FK(\dim M)$ the 
Kostant partition corresponding to the isomorphism class of $M$.
We denote by $M_\theta$ an indecomposable represetation corresponding 
to a raiz $\theta\in R^+$, and we denote by $M_\ka$ a representation 
corresponding to a Kostant partition $\ka\in\FK(\al)$.


Let $\bj:C\to S$ be a closed embedding of a smooth curve $C$
into a surface $S$. Let $\BRN_n(C)$ denote the category of nilpotent
representations of the quiver $\TBA_{n-1}$ in the category of cohernet 
sheaves on $S$ with $0$-dimensional support on the curve~$C$.
Every object $T$ of the category $\BRN_n(C)$ can be decomposed 
as $T=\oplus_{x\in C}T_x$, where $T_x$ is an object concentrated 
at a point $x\in C$. Let 
$$
\Ga, \Ga_x:\BRN_n(C)\to\BRN_n,\qquad \Ga(T)=\Ga(S,T),\quad \Ga_x(T)=\Ga(S,T_x)
$$ 
denote the functor of global sections and of global sections
with support at $x$ respectively.


Let $\ka_x$ be the Kostant partition, corresponding to the isomorphism
class of the object $\Ga_x(T)$ of the category $\BRN_n$. Then the
nontrivial Kostant partitions $\ka_x$, $x\in C$ form a multipartition 
$\mu(T)\in\FM(\al)$ of $\al=\dim\Ga(T)$. We call the objects 
$T$ and $T'$ of $\BRN_n(C)$ {\em equivalent} if $\mu(T)=\mu(T')$.

Thus the set of equivalence classes of objects of the category $\BRN_n(C)$
are in one-to-one correspondence with multipartitions.

\section{The space $K_\al$}

\subsection{Definition and piecification}

We fix a smooth surface $S$ and a smooth curve $C\subset S$.
Let $\SO=S - C$ be the complement, $\bi:\SO\to S$ the open embedding,
and $\bj:C\to S$ the closed embedding. Let $[C]\in H^2(S,\ZZ)$ be the
fundamental class of the marked curve $C\subset S$. We denote by
$d=[C]^2=\deg\CN_{C/S}$ the degree of the normal bundle
and by $g=g(C)$ the genus of the curve $C$.


Let $V$ be a $n$-dimensional vector space. We fix a flag 
$$
\dots\subset\CE_{-1}\subset\CE_0\subset\CE_1\subset\dots\subset V\otimes \bi_*\CO_\SO
$$ 
of rank $n$ vector bundles on the surface $S$ such that 
$$
\rnc{\arraystretch}{1.2}
\begin{array}{ll}
\CE_{p-n}=\CE_p(-C)=\CE_p\otimes\CO_S(-C)\subset \CE_p\qquad\qquad\qquad\qquad &
\text{(periodicity)}\\
c_1(\CE_{p+1})=c_1(\CE_p)+[C]	&\text{(normalization)}
\end{array}
$$
It follows that $\CE_p/\CE_{p-1}=\bj_*\CL_p$, where
$\CL_p$ are line bundles on the curve $C$. Moreover,
the periodicity implies that 
$$
\CL_{p+n}=\CL_p\otimes\bj^*\CO_S(C)=\CL_p\otimes\CN_{C/S}
$$ 

Let $\al_0$ denote the sequence 
formed by $-ch_2(\CE_p)$, where $ch_2=c_1^2/2-c_2$ 
is the second coefficient of the Chern character.

\begin{lm}\label{centralcharge}
We have $\al_0\in Y^{(2)}\otimes\QQ$ and $\langle i',\al_0\rangle=\deg\CL_{i+1}-\deg\CL_i$.
Moreover, we have $\sum_{i\in I}\langle i',\al_0\rangle=d$.
\end{lm}
\begin{pf}
We have
$$
ch(\CE_{p+n})-ch(\CE_p)=ch(\CE_p)(1+[C]+[C]^2/2)-ch(\CE_p)=ch(\CE_p)([C]+d/2[\point]),
$$
hence
$$
ch_2(\CE_{p+n})-ch_2(\CE_p)=ch_1(\CE_p)\cdot[C]+ch_0(\CE_p)\cdot d/2=
pd+c_1(\CE_0)\cdot[C]+nd/2,
$$
hence $\al_0\in Y^{(2)}\otimes\QQ$.

On the other hand by the Riemann-Roch-Grothendick Theorem we have
$$
ch(\bj_*\CL_p)=[C]+(\deg\CL_p-d/2)[\point],
$$
hence 
$$
\deg\CL_p=ch_2(\CE_p)-ch_2(\CE_{p-1})+d/2
$$
hence
$$
\deg\CL_{p+1}-\deg\CL_p=ch_2(\CE_{p+1})+ch_2(\CE_{p-1})-2ch_2(\CE_p)=
\langle p',\al_0\rangle.
$$

Finally, we have
$$
\sum_{i\in I}\langle i',\al_0\rangle=\sum_{p=0}^{n-1}(\deg\CL_{p+1}-\deg\CL_p)=
\deg\CL_n-\deg\CL_0=d.
$$
\end{pf}

Any infinite flag of coherent sheaves on the surface $S$
can be considered as a representation of the quiver $\BA_\infty$.
Given a periodic subflag $E_\bullet\subset\CE_\bullet$, such that
${E_\bullet}_{|\SO}\cong{\CE_\bullet}_{|\SO}$ we denote by
$T_\bullet=\CE_\bullet/E_\bullet$ the quotient representation
of the quiver~$\BA_\infty$. Assume that the support of $T_\bullet$
is $0$-dimensional. Then choosing a trivialization of the normal
bundle $\CN_{C/S}$ in a neighbourhood of $\supp T_\bullet$
we obtain an isomoprhism
$$
T_{p+n}\cong T_n,
$$
hence $T_\bullet$ can be considered as a representation of the cyclic quiver $\TBA_{n-1}$
in the category of coherent sheaves on the curve $C$. Then $\Ga(T_\bullet)$ 
is a representation of $\TBA_{n-1}$ in the category of vector spaces.
It is clear that both $T_\bullet$ and $\Ga(T_\bullet)$ are nilpotent.

\begin{df}
Let $K_\al(\CE_\bullet)$ denote the space of all periodic subflags
$E_\bullet\subset\CE_\bullet$ such that $T_\bullet$ is an
object of the category $\BRN_n(C)$ and $\dim\Ga(T_\bullet)=\alpha$.
\end{df}


We will denote the space $K_\al(\CE_\bullet)$ simply by $K_\al$ for brevity.


For a multipartition $\mu(E_\bullet)\in\FM(\al)$ let $K_\mu\subset K_\al$ 
denote the subspace of all $E_\bullet$ such that the equivalence class
$\mu(T_\bullet)$ of the object $T_\bullet=\CE_\bullet/E_\bullet$ of the
category $\BRN_n(C)$ is equal to $\mu$. This defines a piecification
\begin{equation}\label{strofk}
K_\al=\bigsqcup_{\mu\in\FM(\al)}K_\mu
\end{equation}

\begin{rm}
It is clear that the equivalence class $\mu(T_\bullet)$
doesn't depend on the choice of the trivialization of the
normal bundle $\CN_{C/S}$ involved.
\end{rm}


Let $\ka_x(T)$ denote the isomorphism class of the representation $\Ga_x(T)$. 
Then by definition the multipartition $\mu(T_\bullet)$ is formed by nontrivial 
Kostant partitions $\ka_x(T_\bullet)$. 
Hence we have a map
$$
\sigma:K_\mu\to C^\al_\mu,\qquad
E_\bullet\mapsto\sum_{x\in C}\ka_x(\CE_\bullet/E_\bullet)x.
$$

Let $\mu=\lbr\ka_1,\dots,\ka_m\rbr$. Given an element 
$\sum\ka_rx_r\in C^\al_\mu$ let $F_\mu(\sum\ka_rx_r)$ 
denote the fiber $\sigma^{-1}(\sum\ka_rx_r)\subset K_\mu$.

\begin{lm}\label{factorization}
The map $\sigma$ is a locally trivial fibration.
Moreover, we have an isomorphism
$$
F_\mu(\ka_1x_1+\dots+\ka_mx_m)=F_{\lbr\ka_1\rbr}(\ka_1x_1)\times\dots\times
F_{\lbr\ka_m\rbr}(\ka_mx_m).
$$
\end{lm}
\begin{pf}
Evident.
\end{pf}

Thus the description of the stratum $K_\mu\subset K_\al$
reduces to the description of the space $F_{\lbr\ka\rbr}(\ka x)$
which is called {\em a simple fiber}.

\subsection{Simple fiber}\label{simplefiber}

This subsection is devoted to the proof of the following Theorem.
We fix a point $x\in C$ and a Kostant partition $\ka$.
We denote the simple fiber $F_{\lbr\ka\rbr}(\ka x)$ by $F_\ka$ for brevity.

Recall that 
$$
F_\ka=\{E_\bullet\subset\CE_\bullet\ |\ \supp(\CE_\bullet/E_\bullet)=\{x\}
\text{ and }\ka(\Ga(\CE_\bullet/E_\bullet))=\ka\}
$$

\begin{th}\label{fka}
The space $F_\ka$ is a pseudoaffine space of dimension $||\ka||-\sK(\ka)$.
\end{th}

It is clear that the space $F_\ka$ depends only on the local properties of
the surface $S$ near the point $x$. So, in this subsection we can and will 
replace $S$ by a small neighbourhood of $x$. This allows to fix a trivialization
of the normal bundle $\CN_{C/S}$, giving an isomorphisms 
$$
\CL_{p+n}=\CL_p.
$$

We fix some integer $s$.

Let $\ka_p^q$ ($p\le q$, $s\le q\le s+n-1$) be the coordinates of the
Kostant partition $\ka$ with respect to $R^+_s$. 


For any collection $\ka_p^q$ and subsets $I,J\subset\ZZ$ we define 
$$
\ka_I^J=\sum_{p\in I,\ q\in J}\ka_p^q.
$$

Another possible definition of the collection $\ka_p^q$ is given 
by the following Lemma.

\begin{lm}\label{kafromn}
Assume that $N$ is a $n$-periodic representation of the quiver $\BA_\infty$
such that its isomorphism class in the category of representations of the
quiver $\TBA_{n-1}$ is equal to $\ka$. Let $\hat N_p$ denote the kernel of the 
map $N_p\to N_{s+n}$. Then for all $p\le q\le s+n-1$ we have
$$
\rk(\hat N_p\to \hat N_q)=\ka_{\le p}^{\ge q}.
$$
\end{lm}
\begin{pf}
The Lemma follows from the calculation of the contributions of 
the summands $M_\theta$ ($\theta=(p',q')\in R^+_s$) of $N$ into 
the rank of the above map. It suffices to note that $M_\theta$ 
contributes to the rank of the map iff $p'\le p$ and $q\le q'$.
\end{pf}

Assume that $E_\bullet\in F_\ka$ and let $T_\bullet=\CE_\bullet/E_\bullet$.

Let us denote the torsion sheaf $L^1\bj^*T_s$ by $R$.
Now we will introduce a pair of filtrations on $R$.

The first of them can be constructed quite easily.

\begin{lm}\label{restr}
There are natural isomorphisms
$$
\bj^*T_s=\Coker (T_{s-n}\to T_s),\qquad
L^1\bj^*T_s=\Ker(T_s\to T_{s+n}).
$$
\end{lm}
\begin{pf}
Evident.
\end{pf}

Let us denote 
$$
R_i=\Ker(T_s\to T_{s+i})\qquad(0\le i\le n).
$$
This defines the {\em right} filtration
$$
0=R_0\subset R_1\subset\dots\subset R_{n-1}\subset R_n=R
$$
of the sheaf $R$.

The second one is a little bit more complicated.
First we consider a filtration on the sheaf $T_s$
formed by the sheaves
$$
T_i^s=\CE_i/(E_s\cap\CE_i)=\Im(T_i\to T_s)\subset T_s.
$$
Then we consider 
$$
R^i=L^1\bj^*(\CE_i/(E_s\cap\CE_i))=L^1\bj^*T_i^s\subset L^1\bj^*T_s=R.
$$
This defines the {\em left} filtration
$$
\dots\subset R^{s-2}\subset R^{s-1}\subset R^s=R.
$$
of the sheaf $R$.

\begin{rem}\label{onlyes}
Note that the left filtration $R^\bullet$ of the sheaf $R$ is defined 
by the subsheaf $E_s\subset\CE_s$ only.
\end{rem}

\begin{lm}\label{kafromr}
We have
\begin{equation}\label{rank}
\rk(\Ga(R^p)\to\Ga(R)\to\Ga(R/R_i))=\ka_{\le p}^{\ge s+i},
\qquad(p\le s,\ 0\le i < n).
\end{equation}
\end{lm}
\begin{pf}
Let $N_\bullet=\Ga(T_\bullet)$. Then according to the Lemma~\ref{kafromn}
we have $\ka_{\le p}^{\ge s+i}=\rk(\hat N_p\to\hat N_{s+i})$.
Since $\Ga(\bullet)$ is an exact functor on the category of
torsion sheaves on the curve it follows that
$$
\hat N_p=\Ga\left(\frac{\CE_p\cap E_s(C)}{E_p}\right),\qquad
\hat N_{s+i}=\Ga\left(\frac{\CE_{s+i}\cap E_s(C)}{E_{s+i}}\right).
$$
It is evident that 
$$
\Im\left(\frac{\CE_p\cap E_s(C)}{E_p}\to\frac{\CE_{s+i}\cap E_s(C)}{E_{s+i}}\right)=
\frac{\CE_p\cap E_s(C)}{\CE_p\cap E_{s+i}},
$$
hence 
\begin{equation}\label{r}
\dim\Ga\left(\frac{\CE_p\cap E_s(C)}{\CE_p\cap E_{s+i}}\right)=\ka_{\le p}^{\ge s+i}
\qquad\text{for all $p\le s$, $0\le i<n$.}
\end{equation}
On the other hand, it is easy to show that
$$
R^p=\frac{\CE_p\cap E_s(C)}{\CE_p\cap E_s},\qquad
R_i=\frac{\CE_s\cap E_{s+i}}{E_s}\quad\text{hence}\quad
R/R_i=\frac{\CE_s\cap E_s(C)}{\CE_s\cap E_{s+i}}.
$$
and the image of the map $R^p\to R/R_i$ is equal to
$$
\frac{\CE_p\cap E_s(C)}{\CE_p\cap E_{s+i}}
$$
and the Lemma follows.
\end{pf}

For every $t\le s$ let $X^t_\ka$ denote the space of 
all subsheaves $E\subset\CE_t$ with a filtration
$$
0=R_0\subset R_1\subset\dots\subset R_{n-1}\subset R_n=R=L^1\bj^*(\CE_t/E)
$$
such that the quotient $\CE_t/E$ is concentrated at the point $x$
and the filtration $R_\bullet$ together with the left filtration
$R^p:=L^1\bj^*(\CE_p/(E\cap\CE_p))$ satisfy the condition 
\begin{equation}\label{rank1}
\rk(\Ga(R^p)\to\Ga(R)\to\Ga(R/R_i))=\ka_{\le p}^{\ge s+i},
\qquad(p\le t,\ 0\le i < n).
\end{equation}

We have an obvious map $\pi_s:F_\ka\to X^s_\ka$ sending
$E_\bullet$ to $(E_s,R_\bullet)$, where $R_\bullet$ is
the right filtration of the sheaf $R$.


Thus the problem of description of the space $F_\ka$
reduces to the description of the space $X^s_\ka$
and to the description of the fiber of the map $\pi_s$.

We begin with some notation. 
Assume that $(E,R_\bullet)\in X^t_\ka$ for some $t\le s$.
We denote the intersection $\CE_i\cap E$
by $E^t_i$ and the quotient $\CE_i/E^t_i$ by $T^t_i$. Then we have
a filtration
$$
\dots\subset T^t_{t-2}\subset T^t_{t-1}\subset T^t_t=:T.
$$

\begin{lm}\label{jti}
We have $\bj^*T^t_i=T^t_i/T^t_{i-n}$ for all $i\le t$.
\end{lm}
\begin{pf}
Since $T^t_i=\CE_i/E^t_i$, it follows that $\bj^*T^t_i$ is
isomorphic to the cokernel of the map $T^t_i\to T^t_i$,
induced by the embeddings $\CE_i(-C)=\CE_{i-n}\subset\CE_i$
and $E^t_i(-C)\subset E^t_i$. However, since we have
$$
E^t_i(-C)\subset E^t_{i-n}\subset\CE_{i-n}
$$
it follows that the morphism $T^t_i\to T^t_i$ factors as
the composition of the surjection $T^t_i\to T^t_{i-n}$ 
and the embedding $T^t_{i-n}\to T^t_i$. The Lemma follows.
\end{pf}

\begin{lm}\label{dimtis}
We have 
$$
\dim\Ga(T^t_p)=\sum_{r=0}^\infty\ka_{\le p-rn}^{\ge s}.
$$
\end{lm}
\begin{pf}
Follows from~\ref{jti} and~\refeq{rank1} since 
$\dim\Ga(\bj^*T_i^t)=\dim\Ga(L^1\bj^*T_i^t)=\dim\Ga(R^i)$.
\end{pf}

Consider the restriction of the embedding $\CE_{i-n}\to\CE_{i-1}$
to the curve $C$. Since $\bj^*\CE_{i-1}=\CE_{i-1}/\CE_{i-n-1}$
it follows that the map $\bj^*\CE_{i-n}\to\bj^*\CE_{i-1}$
factors as $\bj^*\CE_{i-n}\to\CL_{i-n}\to\bj^*\CE_{i-1}$.
Thus, we have an embedding $\CL_i=\CL_{i-n}\to\bj^*\CE_{i-1}$.


\begin{lm}\label{kernel}
The kernel of the composition $\CL_i\to \bj^*\CE_{i-1}\to \bj^*T^t_{i-1}$
is equal to $\CL_i(-\sum_{r=1}^\infty\ka_{i-rn}^{\ge s}x)$.
\end{lm}
\begin{pf}
It suffices to show that the image $T$ of the above map
is a torsion sheaf with $\dim\Ga(T)=\ka_{\le i-n}^{\ge s}$. 
Note that the commutative diagram
$$
\begin{CD}
\bj^*\CE_{i-n} @>>> \bj^*\CE_{i-1} \\
@VVV              @VVV         \\
\bj^*T^t_{i-n} @>>> \bj^*T^t_{i-1}
\end{CD}
$$
implies $T=\Im(\bj^*T^t_{i-n} \to \bj^*T^t_{i-1})$.
>From the Lemma~\ref{jti} it follows that 
$$
T=\Im(T^t_{i-n}/T^t_{i-2n} \to T^t_{i-1}/T^t_{i-n-1})=
T^t_{i-n}/T^t_{i-n-1},
$$
hence according to the Lemma~\ref{dimtis} we have 
$$
\dim\Ga(T)=\dim\Ga(T_{i-n}^t)-\dim\Ga(T_{i-n-1}^t)=
\sum_{r=1}^\infty\ka_{\le i-rn}^{\ge s}-\sum_{r=1}^\infty\ka_{\le i-rn-1}^{\ge s}=
\sum_{r=1}^\infty\ka_{i-rn}^{\ge s}.
$$
\end{pf}

\begin{cor}\label{obs}
The composition 
$$
\CL_i\left(-\sum_{r=0}^\infty\ka_{i-rn}^{\ge s}x\right)\to
\CL_i\to \bj^*\CE_{i-1}\to \bj^*T^t_{i-1}
$$
vanishes.
\end{cor}
\begin{pf}
Follows from the above Lemma.
\end{pf}

Consider the map $\varpi_t:X^t_{\ka}\to X^{t-1}_{\ka}$
given by
$$
\varpi_t(E,R_\bullet\subset R^t)=
(E\cap\CE_{t-1},(R^{t-1}\cap R_\bullet)\subset R^{t-1}).
$$

\begin{pr}\label{step}
The map $\varpi_t$ is a locally trivial pseudoaffine fibration with
the fiber $\varpi_t^{-1}(E',R'_\bullet)$ being a pseudoaffine 
space of dimension $\ka_{\le t-1}^{\ge s}$.
\end{pr}

The proof of the Proposition \ref{step} consists of a few steps.

Consider the space $Y$ of all pairs $(E,R'_\bullet)$, where
$E$ is a subsheaf in $\CE_t$ and $R'_\bullet$ is a filtration
in the sheaf $R'=L^1\bj^*(\CE_{t-1}/(E\cap\CE_{t-1}))$ such that
$\CE_t/E$ is concentrated at $x$ and
$\dim\Ga(\CE_t/E)=\sum\limits_{r=0}^\infty\ka_{\le t-rn}^{\ge s}$
and for all $p\le t-1$ the conditions \refeq{rank1} are satisfied.

We have natural maps $\varpi'_t:X^t_\ka\to Y,\ 
(E,R_\bullet)\mapsto (E,R_\bullet\cap R')$, and
$\xi:Y\to X^{t-1}_\ka,\ (E,R'_\bullet)\mapsto (E\cap\CE_{t-1},R'_\bullet)$,
and evidently $\varpi_t=\xi\cdot\varpi'_t$.

We begin with the description of the space $Y$.

\begin{pr}\label{y}
The space $Y$ is a torsor over the vector bundle
$$
H:=\Hom_{X^{t-1}_\ka}(\CL_t(-kx),R'), 
$$
where $k=\sum\limits_{r=0}^\infty\ka_{t-rn}^{\ge s}$.
\end{pr}
\begin{pf}
Let $(E',R'_\bullet)\in X^{t-1}_\ka$ and let $T'$ denote the quotient 
$\CE_{t-1}/E'$. 
If $E$ is a subsheaf in $\CE_t$ such that $E\cap\CE_{t-1}=E'$ and 
$(E,R'_\bullet)\in Y$, then $E/E'$ is a subsheaf in 
$\CE_t/\CE_{t-1}=\bj_*\CL_t$, and 
$$
\dim\Ga\left(\frac{\CE_t/\CE_{t-1}}{E/E'}\right)=
\dim\Ga\left(\frac{\CE_t/E}{\CE_{t-1}/E'}\right)=
\sum\limits_{r=0}^\infty\ka_{\le t-rn}^{\ge s}-
\sum\limits_{r=0}^\infty\ka_{\le t-1-rn}^{\ge s}=k,
$$
hence $E/E'=\bj_*\CL_t(-kx)$. 
Let $\TW$ denote  the quotient $\CL_t/\CL_t(-kx)$,
and let $\TE$ denote the kernel of the composition
$\CE_t\to \bj_*\CL_t\to\TW$. It follows that $Y$ is the space
of all extensions of the projection $\CE_{t-1}\to T'$
to the map $\CE_{t-1}\subset\TE\to T'$. 
It follows from~\cite{mrl} that we have the obstruction
map $\eta:X^{t-1}_\ka\to\Ext^1(\bj_*\CL_t(-kx),T')$ and that
$Y$ is a $\Hom(\bj_*\CL_t(-kx),T')$-torsor over the zero locus 
of the map $\eta$. Note that 
$$
\Hom(\bj_*\CL_t(-kx),T')\cong\Hom(\CL_t(-kx),L^1\bj^*T')=H.
$$
Hence it suffices to prove that $\eta=0$ in our case.

To this end note that the obstruction $\eta(E')$ is equal 
to the Ioneda product of the embedding $\bj_*\CL_t(-kx)\to \bj_*\CL_t$,
of the extension 
\begin{equation}\label{ext}
0 \to \CE_{t-1} \to \CE_t \to \bj_*\CL_t \to 0,
\end{equation}
and of the projection $\CE_{t-1}\to T'$.
On the other hand, we have an isomorphism
$$
\Ext^1(\bj_*\CL_t,\CE_{t-1})\cong \Hom(\CL_t,\bj^*\CE_{t-1})
$$
under which the extension \refeq{ext} corresponds
to the natural embedding $\CL_t\to \bj^*\CE_{t-1}$
(see Lemma \ref{kernel}). Hence the Corollary \ref{obs}
implies the vanishing of the obstruction.
\end{pf}

Assume that $(E,R'_\bullet)\in Y$ and let $E'=E\cap\CE_{t-1}$.
Let $T=\CE_t/E$ and $T'=\CE_{t-1}/E'\subset T$.
It follows from the proof of the Proposition~\ref{y} 
that we have an exact sequence
$$
0 \to T' \to T \to \bj_*(\CL_t/\CL_t(-kx)) \to 0,
$$
hence $R/R'$ is a subsheaf in 
$L^1\bj^*(\bj_*(\CL_t/\CL_t(-kx)))=\CL_t/\CL_t(-kx)$.
On the other hand we have $\dim\Ga(R')=\ka^{\ge s}_{\le t-1}$
by the definition of the space $Y$ and 
\begin{multline*}
\dim\Ga(R)=\dim\Ga(L^1\bj^*T)=\dim\Ga(\bj^*T)=\dim\Ga(T/T^t_{t-n})=\\=
\dim\Ga(T)-\dim\Ga(T^t_{t-n})=\dim\Ga(T)-\dim\Ga(T^{t-1}_{t-n})=\ka_{\le t}^{\ge s}.
\end{multline*}
Hence $\dim\Ga(R/R')=\ka_{\le t}^{\ge s}-\ka^{\ge s}_{\le t-1}=\ka_t^{\ge s}$.
This means that 
$$
R/R'=W:=\CL_t((\ka_t^{\ge s}-k)x)/\CL_t(-kx)=
\CL_t((\ka_t^{[s,s+n-1]}-k)x)/\CL_t(-kx).
$$

Now, if $R_\bullet$ is a filtration on $R$ such that
$(E,R_\bullet)\in X^t_\ka$ then we have 
$R_i/R'_i\subset R/R'=W$ and
$$
\dim\Ga(R_i/R'_i)=\dim\Ga(R_i)-\dim\Ga(R'_i)=
\ka_{\le t}^{[s,s+i-1]}-\ka_{\le t-1}^{[s,s+i-1]}=\ka_t^{[s,s+i-1]},
$$
hence
$$
R_i/R'_i=W_i:=\CL_t((\ka_t^{[s,s+i-1]}-k)x)/\CL_t(-kx).
$$

Thus we have the following standard commutative diagram:
\begin{equation}\label{eps}
\minCDarrowwidth=18pt
\begin{CD}
  @.	0	 	@.   	0	    	@.	0		\\
@.	@VVV	      		@VVV	 		@VVV		\\
  @.   	R'_\bullet 	@.   	R_\bullet 	@.	W_\bullet	\\
@.	@VVV	      		@VVV			@VVV		\\
0 @>>>	R' 		@>>>   	R 		@>>>	W	@>>> 0 	\\
@.	@V{\can}VV		@VVV			@VVV		\\
0 @>>>	R'/R'_\bullet	@>>>	R/R_\bullet	@>>>	W/W_\bullet  @>>> 0	\\	
@.	@VVV	      		@VVV	 		@VVV		\\
  @.	0	 	@.   	0	    	@.	0			
\end{CD}
\end{equation}

\begin{rem}
This is a diagram in the category of representations of the quiver
$$
\nc{\vertex}[1]{\stackrel{#1}{\bullet}}
\nc{\arrow}{\ \longrightarrow\ }
\vertex{0}\arrow\vertex{1}\arrow\dots\arrow\vertex{n-1}\arrow\vertex{n}
$$
in the category of coherent sheaves on the curve $C$. From now and till
the end of the proof of the Proposition the index $\bullet$ indicates an
object of this category. It will be rather important below that indices
start from $0$.
\end{rem}

Let $\TR_\bullet$ denote the kernel of the composition 
$R\to W\to W/W_\bullet$. From the standard technique of \cite{mrl}
it follows that 
$$
X^t_\ka=\left\{f\in\Hom(\TR_\bullet,R'/R'_\bullet)\ \left|\ 
\text{such that the triangle}\
{\def\objectstyle{\textstyle}\def\labelstyle{\scriptstyle}\vcenter{
\xymatrix @-.95pc {
{R'} \ar[r] \ar[d]_{{\can}} & {\TR_\bullet} \ar[dl]^{f} \\
{R'/R'_\bullet}
}}}
\text{commutes}
\right.\right\}
$$

In other words, we have the following cartesian square
$$
\begin{CD}
X^t_\ka		@>>>	\Hom_Y(\TR_\bullet,R'/R'_\bullet)	\\
@VVV			@VVV					\\
Y	@>1\times\can>>	Y\times\Hom_{X^{t-1}_\ka}(R',R'/R'_\bullet)
\end{CD}
$$

Consider the map $\eps:Y\to\Ext^1(W,R')$ given
by the middle row of \refeq{eps} and the projection
$\Ext^1(W,R')\to\Ext^1(W_\bullet,R')$
induced by the embedding $W_\bullet\to W$. Note that applying
the functor $\Hom(\bullet,R')$ to the sequences
\begin{equation}\label{s}
0\to\CL_t(-kx)\to\CL_t((\ka_t^{[s,s+n-1]}-k)x)\to W\to0,
\end{equation}
\begin{equation}\label{sb}
0\to\CL_t(-kx)\to\CL_t((\ka_t^{[s,s+\bullet-1]}-k)x)\to W_\bullet\to0
\end{equation}
we obtain the morphisms of $X^{t-1}_\ka$-spaces
$$
\Hom(\CL_t(-kx),R')\to\Ext^1(W,R')\quad\text{and}\quad
\Hom(\CL_t(-kx),R')\to\Ext^1(W_\bullet,R')
$$
which define a natural fiberwise (over $X^{t-1}_\ka$) action 
of the vector bundle $H$ on $\Ext^1(W,R')$ and $\Ext^1(W_\bullet,R')$. 

\begin{lm}
The map $\teps:Y\stackrel{\eps}{\longrightarrow}\Ext^1(W,R')\longrightarrow\Ext^1(W_\bullet,R')$
commutes with the action of $H$.
\end{lm}
\begin{pf}
Let $\TW$ denote the quotient $\CL_t/\CL_t(-kx)$. 
Recall that we have the following commutative diagrams
$$
\begin{CD}
0 @>>> \TE 	       @>>> \CE_t      @>>> \bj_*\TW @>>> 0	\\
@.	@VVV		    @VVV		@|		\\
0 @>>> \bj_*\CL_t(-kx) @>>> \bj_*\CL_t @>>> \bj_*\TW @>>> 0
\end{CD}
$$
and
$$
\begin{CD}
0 @>>> \CL_t(-kx) @>>> \CL_t((\ka_t^{[s,s+\bullet-1]}-k)x) @>>> W_\bullet	@>>> 0	\\
@.	@|		@VVV				     	@VVV			\\
0 @>>> \CL_t(-kx) @>>> \CL_t((\ka_t^{[s,s+n-1]}-k)x) 	@>>> 	W   		@>>> 0	\\
@.	@|		@VVV				     	@VVV			\\
0 @>>> \CL_t(-kx) @>>> \CL_t 				@>>> 	\TW 		@>>> 0	
\end{CD}
$$
Applying the functors $\Hom(\bullet,T')$ and $\Hom(\bullet,R')$
we get the following commutative diagram
$$
\xymatrix@C-15pt{
Y \ar[d] & \Hom(\bj_*\CL_t(-kx),T') \ar[dl] \ar[d] \ar@{=}[r] &
\Hom(\CL_t(-kx),R') \ar[d] \ar[dr] \ar[r] & \Ext^1(W_\bullet,R')\\
\Hom(\TE,T') \ar[r] & \Ext^1(\bj_*\TW,T') \ar@{=}[r] & 
\Ext^1(\TW,R') \ar[r] & \Ext^1(W,R') \ar[u]
}
$$
and the Lemma follows.
\end{pf}


\begin{cor}\label{comm}
We can choose a local over $X^{t-1}_\ka$ trivialization 
$\phi:Y\stackrel{\sim}{\to} H$ such that 
the following diagram
$$
\xymatrix{
Y \ar[r]^{\phi} \ar[d]_{\teps} & H \ar[dl] \\
\Ext^1(W_\bullet,R')
}
$$
commutes.
\end{cor}

Thus for every point $(E,R'_\bullet)\in Y$ we have a homomorphism 
$\phi_E:\CL_t(-kx)\to R'$.
\begin{lm}\label{psi}
The homomorphism $\phi_E$ can be extended $($locally over $X^{t-1}_\ka)$
to a morphism of complexes
$$
\begin{CD}
0  @>>>	\CL_t(-kx) @>>>	\CL_t((\ka_t^{[s,s+\bullet-1]}-k)x) @>>> W_\bullet @>>> 0 \\
@.	@V{\phi_E}VV	@V{\psi_E}VV				 @|		  \\
0  @>>>	R'	   @>>>	\TR_\bullet			    @>>> W_\bullet @>>> 0
\end{CD}
$$
\end{lm}
\begin{pf}
The claim of the Corollary \ref{comm} reformulated in terms of the derived 
category says that the square in the following diagram
$$
\begin{CD}
W_\bullet[-1]	@>>> \CL_t(-kx) @>>> \CL_t((\ka_t^{[s,s+\bullet-1]}-k)x) 
@>>> W_\bullet\\
@|				@V{\phi_E}VV 	@.	@|\\
W_\bullet[-1]	@>{\teps_E}>>	R'	     @>>> \TR_\bullet @>>> W_\bullet
\end{CD}
$$
commutes. Hence it can be extended (locally) to a morphism of triangles,
and the space of extensions is a torsor over $\Hom(W_\bullet,\TR_\bullet)$.
Hence the obstruction to the local over $X^{t-1}_\ka$ extension lies
in $R^1\xi_*(\Hom(W_\bullet,\TR_\bullet))$. But the map $\xi$ is an affine
morphism, hence the obstruction vanishes.
\end{pf}

Thus we have a map $\Hom_Y(\TR_\bullet,R'/R'_\bullet)\stackrel{\psi}{\to}
\Hom(\CL_t((\ka_t^{[s,s+\bullet-1]}-k)x),R'/R'_\bullet)$
induced by the morphisms $\psi_E$. Consider the following diagram
$$
\begin{CD}
\Hom_Y(\TR_\bullet,R'/R'_\bullet) 				@>{\psi}>> 
\Hom(\CL_t((\ka_t^{[s,s+\bullet-1]}-k)x),R'/R'_\bullet)			\\
@VVV			@VVV							\\
Y\times\Hom(R',R'/R'_\bullet) 			@>{\phi}>> 
\Hom(\CL_t(-kx),R'/R'_\bullet)
\end{CD}
$$

\begin{lm}
The above square is cartesian.
\end{lm}
\begin{pf}
The Lemma \ref{psi} implies that the above square commutes, so it suffices to 
note that the fibers of the map $\Hom_Y(\TR_\bullet,R'/R'_\bullet) \to 
Y\times\Hom(R',R'/R'_\bullet)$ are equal to $\Hom(W_\bullet,R'/R'_\bullet)$
and that they map isomorphically to the fibers of the map
$\Hom(\CL_t((\ka_t^{[s,s+\bullet-1]}-k)x),R'/R'_\bullet)\to
\Hom(\CL_t(-kx),R'/R'_\bullet)$.
\end{pf}

\begin{cor}
We have a cartesian square
$$
\begin{CD}
X^t_\ka	@>>>		\Hom(\CL_t((\ka_t^{[s,s+\bullet-1]}-k)x),R'/R'_\bullet) \\
@VVV			@VVV							\\
Y 	@>{\phi}>> 	\Hom(\CL_t(-kx),R'/R'_\bullet)
\end{CD}
$$
\end{cor}

We will need the following Lemma.
\begin{lm}\label{homandext}
$(i)$ $\Hom(\CL_t(-kx),R'_\bullet)=\Ext^1(\CL_t(-kx),R'_\bullet)=0$.

\medskip
\noindent
$(ii)$ The projection $R'\to R'/R'_\bullet$ induces an isomorphism
$$
\Hom(\CL_t(-kx),R')\stackrel\sim\to\Hom(\CL_t(-kx),R'/R'_\bullet).
$$

\noindent
$(iii)$ 
$\Ext^1(\CL_t((\ka_t^{[s,s+\bullet-1]}-k)x),R')=
\Ext^1(\CL_t((\ka_t^{[s,s+\bullet-1]}-k)x),R'/R'_\bullet)=0,$
$$
\dim\Hom(\CL_t((\ka_t^{[s,s+\bullet-1]}-k)x),R'/R'_\bullet)=\ka_{\le t-1}^{\ge s}.
$$
\end{lm}
\begin{pf}
Easy.
\end{pf}

Now we can finish the proof of the Proposition \ref{step}. Just note that 
from the definition of $\phi$ and from the Lemma~\ref{homandext} $(ii)$ 
it follows that $\phi$ is a local (over $X^{t-1}_\ka$) isomorphism, 
hence $X^t_\ka$ is locally (over $X^{t-1}_\ka$) isomorphic to the 
vector bundle $\Hom(\CL_t((\ka_t^{[s,s+\bullet-1]}-k)x),R'/R'_\bullet)$,
which according to the Lemma~\ref{homandext} $(iii)$ 
is $\ka_{\le t-1}^{\ge s}$-dimensional.

\begin{rem}\label{allexts}
It follows also that the natural map 
$$
\varpi_t^{-1}(E',R'_\bullet) \to \Ext^1(W_1,R'_1)
$$
is surjective. Indeed, the map
$$
\Hom(\CL_t((\ka_t^{[s,s+\bullet-1]}-k)x),R'/R'_\bullet) \to
\Ext^1(\CL_t((\ka_t^{[s,s+\bullet-1]}-k)x),R'_\bullet)
$$
is surjective by the Lemma~\ref{homandext} $(iii)$.
On the other hand
$$
\Ext^1(\CL_t((\ka_t^{[s,s+\bullet-1]}-k)x),R'_\bullet)\cong\Ext^1(W_\bullet,R'_\bullet)
$$
by the Lemma~\ref{homandext} $(i)$. So it remains to note that the projection
$$
\Ext^1(W_\bullet,R'_\bullet)\to\Ext^1(W_1,R'_1)
$$
is surjective and that $\varpi_t^{-1}(E',R'_\bullet)\cong
\Hom(\CL_t((\ka_t^{[s,s+\bullet-1]}-k)x),R'/R'_\bullet)$.
\end{rem}

Now we can describe the space $X^s_\ka$.

\begin{pr}\label{xs}
The space $X^s_\ka$ is a pseudoaffine space of dimension 
$$
\dim X^s_\ka=\sum_{p\le s}(s-p)\ka_p^{\ge s}=\sum_{p\le s\le q\le s+n-1}(s-p)\ka_p^q.
$$
\end{pr}
\begin{pf}
Follows by induction from the Proposition~\ref{step}.
\end{pf}

Recall that our goal is to describe the space $F_\ka$.
Since we have the map $\pi_s:F_\ka\to X^s_\ka$
and the description of the space $X^s_\ka$ is given by the
Proposition \ref{xs}, it remains to describe the fiber 
of the map $\pi_s:F_\ka\to X^s_\ka$.

We will need the following Lemma.

\begin{lm}\label{dims}
If $E_\bullet\in F_\ka$ then for all $p\le q\le s+n-1$ we have
$$
\dim\Ga\left(\frac{E_{q+1}\cap\CE_p}{E_q\cap\CE_p}\right)=\ka_{\le p}^q.
$$
\end{lm}
\begin{pf}
Follows immediately from \refeq{r}.
\end{pf}

\begin{pr}\label{pis}
The map $\pi_s$ is a pseudoaffine fibration of dimension
$$
\sum_{s\le p < q\le s+n-1}\ka_{\le p}^q=
\sum_{s< p\le q\le s+n-1}(q-p)\ka_p^q+
\sum_{p\le s\le q\le s+n-1}(q-s)\ka_p^q.
$$
\end{pr}
\begin{pf}
Assume that $(E,R_\bullet)\in X^s_\ka$ and 
let $E_\bullet\in\pi_s^{-1}(E,R_\bullet)$.
Then we have a diagram
$$
\begin{array}{ccccccccc}
E_s & \subset & E_{s+1} & \subset & \dots & \subset & E_{s+n-1} & \subset & E_{s+n} \\
\cap && \cap &&&& \cap && \cap \\
\CE_s & \subset & \CE_{s+1} & \subset & \dots & \subset & \CE_{s+n-1} & \subset & 
\CE_{s+n}
\end{array}
$$
But $E_s=E$ and $E_{s+n}=E_s(C)=E(C)$, hence every $E_{s+i}$ is a subsheaf in 
$E(C)\cap\CE_{s+i}$. In other words, we have the following diagram
$$
\begin{array}{ccccccccc}
E & \subset & E_{s+1} & \subset & \dots & \subset & E_{s+n-1} & \subset & E(C) \\
\cap && \cap &&&& \cap && || \\
E(C)\cap\CE_s & \subset & E(C)\cap\CE_{s+1} & \subset & \dots & \subset & 
E(C)\cap\CE_{s+n-1} & \subset & E(C)
\end{array}
$$

Consider the following sheaves:
$$
\TE_i:=\frac{E(C)\cap\CE_{s+i}}E,\qquad
E'_i:=\frac{E_{s+i}}E
$$
(note that $\TE_0=R$). We have a diagram
$$
\begin{array}{ccccccccc}
0 & \subset & E'_1 & \subset & \dots & \subset & E'_{n-1} & \subset & \TE_n \\
\cap && \cap &&&& \cap && || \\
\TE_0 & \subset & \TE_1 & \subset & \dots & \subset & 
\TE_{n-1} & \subset & \TE_n
\end{array}
$$
and it is easy to show that for all $i\le j$ we have
$$
\frac{\TE_i\cap E'_{j+1}}{\TE_i\cap E'_j}=
\frac{\CE_{s+i}\cap E_{s+j+1}}{\CE_{s+i}\cap E_{s+j}},
$$
hence according to \ref{dims} we have
\begin{equation}\label{ij}
\dim\Ga\left(\frac{\TE_i\cap E'_{j+1}}{\TE_i\cap E'_j}\right)=\ka_{\le s+i}^{s+j}.
\end{equation}

Note also that
\begin{equation}\label{cap}
E'_i\cap\TE_0=\frac{\CE_s\cap E_{s+i}}E=R_i\subset R_n=\frac{\CE_s\cap E(1)}E=\TE_0
\end{equation}
(compare with the proof of the Lemma~\ref{kafromr}). Thus, starting from $E_\bullet$
we obtained a flag of subsheaves in the flag $\TE_\bullet$ such that
\refeq{ij} and \refeq{cap} are satisfied. On the other hand it is easy to show 
that if $(E'_\bullet\subset\TE_\bullet)$ is such a flag, then 
$E_{s+i}=\Ker(E(C)\cap\CE_{s+i}\to\TE_i\to \TE_i/E'_i)$
gives a periodic flag $E_\bullet\in\pi_s^{-1}(E,R_\bullet)$.

Thus we have to describe the space of all subflags 
$E'_\bullet\subset\TE_\bullet$ such that \refeq{ij} and \refeq{cap}
are satisfied. Note that all the quotients
$$
\frac{\TE_{i+1}}{\TE_i}\subset\frac{\CE_{s+i+1}}{\CE_{s+i}}
$$
are invertible sheaves on $C$, and that all the intersections
$E'_i\cap\TE_0$ are fixed (according to \refeq{cap}), hence
we can apply the standard technique of \cite{mrl}.
It follows that $\pi_s^{-1}(E,R_\bullet)$ is a pseudoaffine 
space of dimension
\begin{multline*}
\sum_{0\le i<j\le n-1}\dim\Ga\left(\frac{\TE_i\cap E'_{j+1}}{\TE_i\cap E'_j}\right)=
\sum_{0\le i<j\le n-1}\ka_{\le s+i}^{s+j}=\\=
\sum_{s\le p < q\le s+n-1}\ka_{\le p}^q=
\sum_{s< p\le q\le s+n-1}(q-p)\ka_p^q+\sum_{p\le s\le q\le s+n-1}(q-s)\ka_p^q.
\end{multline*}

\end{pf}

Now the Theorem \ref{fka} follows from the Propositions~\ref{xs} and \ref{pis}.

\subsection{Topdimensional components of $K_\al$}

Applying Theorem~\ref{fka}, Lemma~\ref{factorization}
and the definition of the space $C^\al_\mu$ we get the following
Proposition.

\begin{pr}
Let $\mu=\lbr\ka_1,\dots,\ka_m\rbr\in\FM(\al)$. 
The stratum $K_\mu\subset K_\al$ is a smooth variety and
$$
\dim K_\mu=|\al|+(1-\sK(\ka_1))+\dots+(1-\sK(\ka_m)).
$$
Therefore $\dim K_\al=|\al|$ and any topdimensional
irreducible component of the space $K_\al$ coincides with the
closure $\OK_A$ of the stratum $K_A\subset K_\al$, corresponding to
a Kostant partition $A\in\FK(\al)$, considered as a simple 
multipartition.
\end{pr}
\begin{pf}
The first part of the Proposition is evident.

Since $\sK(\ka)\ge 1$ for any Kostant partition $\ka$ it follows
that for any stratum $K_\mu$ we have $\dim K_\mu\le|\al|$, and 
equality is achieved iff $\sK(\ka_1)=\dots=\sK(\ka_m)=1$, that
is iff $\mu$ is a simple multipartition of $\al$. So
it remains to note that a simple multipartition is nothing 
but a Kostant partition.
\end{pf}

\section{The space $\FE_\al^i$}

\subsection{Definition}


Let $\FE_\al^i\subset K_\al\times K_{\alp}$ denote
the closed subspace formed by pairs 
$(E_\bullet,E'_\bullet)\in K_\al\times K_{\alp}$
such that $E'_\bullet\subset E_\bullet$.

The embedding $\FE_\al^i\subset K_\al\times K_{\alp}$ induces
the projections 
$$
\bp:\FE_\al^i\to K_\al\quad\text{and}\quad
\bq:\FE_\al^i\to K_{\alp}.
$$
If $(E_\bullet,E'_\bullet)\in\FE_\al^i$ then 
it is clear that $E_\bullet/E'_\bullet\cong M_i\otimes\CO_x$
as an object of category $\BRN_n(C)$ (here $x\in C$ and
$\CO_x$ stands for the structure sheaf of the point $x$).

Let $\br:\FE_\al^i\to C$ denote the map sending $(E_\bullet,E'_\bullet)$
to the point $x$.

Let $E'_\bullet\in K_{\alp}$ and let $T'_\bullet=\CE_\bullet/E'_\bullet$.

\begin{lm}\label{fiberqr}
The fiber of the map $\bq\times\br:\FE_\al^i\to K_{\alp}\times C$
over the point $(E'_\bullet,x)$ is isomorphic to the projective space
$$
\PP(\Hom(M_i\otimes\CO_x,T'_\bullet))=\PP(\Hom(\CO_x,\Ker(T'_i\to T'_{i+1}))).
$$
\end{lm}
\begin{pf}
Clear.
\end{pf}

Consider a piecification
$$
K_{\alp}\times C=\bigsqcup_{|\ka|\le\alp}Z_{\alp}^\ka,
$$
where $Z_{\alp}^\ka$ denotes the space of pairs $(E'_\bullet,x)$
such that $\ka(\Ga_x(T'_\bullet))=\ka$.


\begin{lm}\label{ztoc}
The projection $Z_{\alp}^\ka\to C$ is a locally trivial fibration.
\end{lm}
\begin{pf}
Easy.
\end{pf}

\begin{lm}\label{zfactor}
We have the following factorization property:
$$
Z_{\alp}^\ka\cong Z_{\al-\be}^0\times_C Z_{\bep}^\ka,
$$
where $\bep=|\ka|$. Moreover, if\/
$W_{\alp}^\ka=(\bq\times\br)^{-1}(Z_{\alp}^\ka)$ then
we have a commutative diagram in which all squares are cartesian
$$
\begin{CD}
\FE_\be^i @<<< W_{\bep}^\ka 	@<<< W_{\alp}^\ka @= 
W_{\alp}^\ka @>>> \FE_\al^i \\
@VVV	@VVV		@VVV		@VVV			@VVV	\\
K_{\bep}\times C @<<< Z_{\bep}^\ka @<<< Z_{\al-\be}^0\times_C Z_{\bep}^\ka @= 
Z_{\alp}^\ka 	@>>> K_{\alp}\times C
\end{CD}
$$
\end{lm}
\begin{pf}
Evident.
\end{pf}

The above Lemma reduces the description of the fibres of the map 
$\bq\times\br:\FE_\al^i\to K_{\alp}\times C$ to the description
of the fibers of the map $W_{\bep}^\ka\to Z_{\bep}^\ka$, where
$|\ka|=\bep$. It is clear that $Z_{\bep}^\ka\to C$ is a locally
trivial $F_\ka$-fibration. We fix some point $x\in C$ and consider
$F_\ka$ as the fiber of $Z_{\bep}^\ka$ over the point $x$. Let
$\TF_\ka=(\bq\times\br)^{-1}(F_\ka)\subset W_{\bep}^\ka$.

Finally, let $Z_{\bep}^\ka(r)\subset Z_{\bep}^\ka$ be the
subspace of points with the dimension of the fiber of the map
$W_{\bep}^\ka\to Z_{\bep}^\ka$ greater than or equal to $r$,
and let $F_\ka(r)\subset F_\ka$ be the subspace of points with 
the dimension of the fiber of the map $\TF_\ka\to F_\ka$ 
greater than or equal to $r$.

\begin{lm}\label{zkar}
The map $Z_{\bep}^\ka(r)\to C$ is a locally trivial $F_\ka(r)$-fibration.
\end{lm}
\begin{pf}
Clear.
\end{pf}

\begin{pr}\label{fkar}
The subspace $F_\ka(r)\subset F_\ka$ is empty if $\sK(\ka)<r+1$
and has codimension at least $r$ if $\sK(\ka)\ge r+1$.
\end{pr}

The proof of the Proposition will be given in the next subsection.
Now we will deduce from it the following Theorem.

\begin{th}\label{ssq}
The map $\bq:\FE_\al^i\to K_{\alp}$ is semismall.
\end{th}
\begin{pf}
It follows from the Proposition \ref{fkar} and from the Lemma \ref{zkar}
that the subspace $Z_{\bep}^\ka(r)\subset Z_{\bep}^\ka$ is empty if 
$\sK(\ka)<r+1$ and has codimension at least~$r$ if $\sK(\ka)\ge r+1$.
>From the Lemma \ref{ztoc} and \ref{zfactor} it follows that the map 
$Z_{\alp}^\ka\to Z_{\bep}^\ka$ is a locally trivial fibration, hence
the subspace $Z_{\alp}^\ka(r)\subset Z_{\alp}^\ka$ is empty if 
$\sK(\ka)<r+1$ and has codimension at least $r$ if $\sK(\ka)\ge r+1$.

Now let $X(r)\subset K_{\alp}\times C$ be the subspace of points 
with the dimension of the fiber of the map $\FE_\al^i\to K_{\alp}\times C$ 
greater than or equal to $r$. Then 
$$
X(r)=\bigsqcup_{|\ka|\le\alp}(X(r)\cap Z_{\alp}^\ka)=
\bigsqcup_{|\ka|\le\alp}Z_{\alp}^\ka(r).
$$
Hence
$$
\codim_{K_{\alp}\times C}X(r)=\min\codim_{K_{\alp}\times C}Z_{\alp}^\ka(r).
$$
But
$$
\codim_{K_{\alp}\times C}Z_{\alp}^\ka(r)=
\codim_{Z_{\alp}^\ka}Z_{\alp}^\ka(r)+
\codim_{K_{\alp}\times C}Z_{\alp}^\ka,
$$
and if $Z_{\alp}^\ka(r)$ is not empty then 
$$
\codim_{K_{\alp}\times C}Z_{\alp}^\ka=\sK(\ka)\ge r+1\quad\text{and}\quad
\codim_{Z_{\alp}^\ka}Z_{\alp}^\ka(r)\ge r,
$$
hence $\codim X(r)\ge 2r+1$.

So it remains to note that the fiber of the map
$\bq:\FE_\al^i\to K_{\alp}$ over the point $E'_\bullet$ 
is equal to the disjoint union of the fibers of the map 
$\bq\times\br$ over the finite number of points $(E'_\bullet,x)$,
hence we have $K_{\alp}(r)\subset p_1(X(r))$, where $K_{\alp}(r)$
is the subspace of points with the dimension of the fiber of the map 
$\FE_\al^i\to K_{\alp}$ greater than or equal to $r$, and 
$p_1:K_{\alp}\times C\to K_{\alp}$ is the projection.
Hence, $\codim K_{\alp}(r)\ge\codim X(r)-1\ge 2r$.
\end{pf}

\subsection{Space $F_\ka(r)$}

This subsection is devoted to the proof of the Proposition~\ref{fkar}.
Here we use the notation of the section~\ref{simplefiber}.
Recall that there we fixed an arbitrary integer $s$.
Let us take $s=i \bmod n$. Then we have 
$$
\Ker(T_i\to T_{i+1})=R_1
$$
according to the definition of $R_1$.
Hence the subspace $F_{\ka}(r)\subset F_\ka$ is given by the condition
$$
\dim\Hom(\CO_x,R_1)\ge r+1.
$$

Recall also that we have a locally trivial fibration $F_\ka\to X^s_\ka$
and a sequence of locally trivial fibrations
$$
X^s_\ka \stackrel{\varpi_s}{\longrightarrow} 
X^{s-1}_\ka \stackrel{\varpi_{s-1}}{\longrightarrow} \dots
$$
Let $X^t_\ka(r)\subset X^t_\ka$ be the subspace of points
$(E,R_\bullet)$ such that $\dim\Hom(\CO_x,R_1)=r+1$.
Let $N_t=\#\{p\le t\ |\ \ka_p^s>0\}$. Then it is clear that $N_s\le \sK(\ka)$.
Hence it suffices to prove the following Proposition.

\begin{pr}
The subspace $X^t_\ka(r)\subset X^t_\ka$ is empty if 
$N_t<r+1$ and has codimension at least $r$ if $N_t\ge r+1$.
\end{pr}
\begin{pf}
We apply the induction in $t$. 

The base of the induction is evident: if $t\ll s$ we have $N_t=0$
and $X^t_\ka(r)$ is empty for any $r$ (because $R_1=0$).

So assume that the induction hypothesis for $t-1$ is true.
We have two cases: $\ka^s_t=0$ and $\ka^s_t>0$.

If $\ka^s_t=0$ then it is evident that $N_t=N_{t-1}$
and that $R_1=R'_1$, hence $X^t_\ka(r)=\varpi_t^{-1}(X^{t-1}_\ka(r))$
and the induction hypothesis for $t-1$ and $t$ are equivalent.

If $\ka^s_t>0$ then $N_t=N_{t-1}+1$ and we have the following exact sequence
\begin{equation}\label{ext1}
0 \to R'_1\to R_1 \to W_1 \to 0
\end{equation}
and $W_1\cong\CL_t/\CL_t(-\ka^s_tx)$. It follows that
$$
X^t_\ka(r)\subset \varpi_t^{-1}(X^{t-1}_\ka(r-1))\bigcup 
\varpi_t^{-1}(X^{t-1}_\ka(r)).
$$
Hence it suffices to check that 
$$
\codim_{\varpi_t^{-1}(X^{t-1}_\ka(r-1))}
(X^t_\ka(r)\cap \varpi_t^{-1}(X^{t-1}_\ka(r-1)))\ge 1.
$$
The latter condition will be satisfied if we prove
that for any point of $X^{t-1}_\ka(r-1)$ a generic point
of the fiber of the map $\varpi_t$ over this point belongs 
to $X^t_\ka(r-1)$.

But this is true, because we have $\dim\Hom(\CO_x,R_1)=\dim\Hom(\CO_x,R'_1)$
for a generic extension \refeq{ext1} and according to the Remark~\ref{allexts} 
all types of extensions \refeq{ext1} are realized
in the fiber of the map $\varpi_t$ over any point of $X^{t-1}_\ka$. 
\end{pf}

\subsection{Irreducible components of $\FE_\al^i$}

We will need a description of irreducible components of $\FE_\al^i$
which are dominant over topdimensional components of the spaces
$K_\al$ and $K_{\alp}$.

We begin with the case of $K_{\alp}$. Let $A'\in\FK(\alp)$ be
a Kostant partition. We will consider $A'$ as a simple multipartition.
Let $K_{A'}$ be the corresponding stratum of $K_{\alp}$ and 
let $\OK_{A'}$ be the corresponding $(|\al|+1)$-dimensional
irreducible component. 

\begin{pr}\label{qdom}
If $A'=\lbr\theta'_1,\dots,\theta'_m\rbr$ then the 
components of $\FE_\al^i$ dominant over $\OK_{A'}$ are in one-to-one 
correspondence with elements $\theta'_r$ of the partition $A'$ 
ending at $i$. The projection $\bp:\FE_\al^i\to K_\al$
sends the generic point of the component of $\FE_\al^i$, corresponding to
the element $\theta'_r$ of the partition $A'$, to the generic
point of the component $\OK_A$ of $K_\al$ with 
$A=\lbr\theta'_1,\dots,\theta'_r\frown i,\dots,\theta'_m\rbr$.
\end{pr}
\begin{pf}
Let $E'_\bullet$ be a generic 
point of $K_{A'}$ and let $\sigma(E'_\bullet)=\theta'_1x_1+\dots+\theta'_mx_m$.
It follows evidently from the Lemma~\ref{fiberqr} that the fiber of the map $\bq\times\br$ over
the point $(E'_\bullet,x)$ is a point, if $x=x_r$ and $\theta'_r\in \sE_i$ (that is if $M_i$ is a 
subrepresentation in $M_{\theta'_r}=\Ga((T'_x)_\bullet)$), and is empty otherwise. Hence, the fiber 
of $\bq$ over the point $E'_\bullet$ is finite and points in the fiber are in a
bijection
with elements of the partition $A'$ ending at $i$.

Let us take a point $(E_\bullet\supset E'_\bullet)\in\FE_\al^i$
in the fiber of $\bq$ over a generic point 
$E'_\bullet\in K_{A'}\subset\OK_{A'}$, corresponding 
to an element $\theta'_r$ of the Kostant partition $A'$.
Then it is clear that for all $x\ne x_r$ we have
$(T_x)_\bullet=(T'_x)_\bullet$ and for $x=x_r$ we have
$$
0 \to M_i\otimes\CO_x \to (T'_x)_\bullet \to (T_x)_\bullet \to 0,
$$
hence $\Ga((T_x)_\bullet)=\Ga((T'_x)_\bullet)/M_i=M_{\theta'_r}/M_i=M_{\theta'_r\frown i}$,
hence $E_\bullet\in K_A$ with $A=\lbr\theta'_1,\dots,\theta'_r\frown i,\dots,\theta'_m\rbr$
and the Proposition follows.
\end{pf}

Let $\FE_{A'}(\theta')$ denote the 
component of $\FE_\al^i$ dominant over $\OK_{A'}$, 
corresponding to the element 
$\theta'$ of the Kostant partition $A'$.

\begin{lm}\label{qdegree}
The map $\bq:\FE_{A'}(\theta')\to\OK_{A'}$ is a generically
finite map of the degree equal to $m(\theta',A')$.
\end{lm}
\begin{pf}
Evident.
\end{pf}

The description of the components of $\FE_\al^i$ which are dominant
over $|\al|$-dimensional components of $K_\al$ is more difficult.
Let $A\in\FK(\al)$ be a Kostant partition. We consider $A$ as a
simple multipartition. Let $K_A$ be the corresponding stratum of 
$K_\al$ and let $\OK_A$ be the corresponding irreducible component.

Let $E_\bullet\in K_A$ be a generic point.

\begin{lm}
The fiber of the map $(\bp\times\br):\FE_\al^i\to K_\al\times C$
over the point $(E_\bullet,x)$ is isomorphic to the projective
space
$$
\PP(\Hom(E_\bullet,M_i\otimes\CO_x))=\PP(\Hom(E_i/E_{i-1},\CO_x)).
$$
\end{lm}
\begin{pf}
Easy.
\end{pf}

Assume that $A=\lbr\theta_1,\dots,\theta_m\rbr$ and
$\sigma(E_\bullet)=\lbr\theta_1\rbr x_1+\dots+\lbr\theta_m\rbr x_m$. 

\begin{lm}
If $x=x_r$ and the raiz $\theta_r\in \sE_{i-1}$ then we have
$$
\PP(\Hom(E_i/E_{i-1},\CO_x))=\PP^1
$$ 
and otherwise $\PP(\Hom(E_i/E_{i-1},\CO_x))$ is a point.
\end{lm}
\begin{pf}
>From the commutative diagram
$$
\begin{CD}
	@.	0		@.	0	\\
@.		@VVV			@VVV	\\
0 	@>>>	E_{i-1}		@>>>	E_i	@>>>	E_i/E_{i-1}	@>>> 0\\
@.		@VVV			@VVV		@VVV		\\
0 	@>>>	\CE_{i-1}	@>>>	\CE_i	@>>>	\bj_*\CL_i	@>>> 0\\
@.		@VVV			@VVV	\\
 	@.	T_{i-1}		@>>>	T_i	\\
@.		@VVV			@VVV	\\
	@.	0		@.	0	
\end{CD}
$$
it follows that $E_i/E_{i-1}$ is a direct sum of a line bundle on the curve $C$
and of the torsion sheaf $\Ker(T_{i-1}\to T_i)$. Hence
$$
\dim\Hom(E_i/E_{i-1},\CO_x)=1+\dim\Hom(\Ker(T_{i-1}\to T_i),\CO_x).
$$
So it remains to note that the sheaf $\Ker(T_{i-1}\to T_i)$ has
a nontrivial component 
at the point $x$ only if $x=x_r$ and $\theta_r\in \sE_{i-1}$,
and that in this case the component is isomorphic to $\CO_x$.
\end{pf}

\begin{cor}
The fiber $\bp^{-1}(E_\bullet)$ is a reducible curve
$$
\BC=C_0\cup\left(\bigcup_{\{r\ |\ \text{$\theta_r\in \sE_{i-1}$}\}}C_r\right),
$$
$($see the Figure $1)$. All vertical components of the curve $\BC$ are genus $0$ curves.
The horizontal component $C_0$ of $\BC$ maps isomorphically to $C$ under
the map $\br$ $($hence $g(C_0)=g(C))$, while the vertical components $C_r$ 
are contracted by $\br$ to the points $x_r\in C$. 
\end{cor}

\begin{figure}[h]
\begin{picture}(340,150)(0,0)
\linethickness{0.5mm}
\put(50,65){\line(1,0){240}}
\put(80,40){\line(0,1){100}}
\put(170,10){\line(0,1){100}}
\put(210,30){\line(0,1){100}}
\put(60,65){\circle*{4}}
\put(80,65){\circle*{4}}
\put(120,65){\circle*{4}}
\put(145,65){\circle*{4}}
\put(170,65){\circle*{4}}
\put(210,65){\circle*{4}}
\put(255,65){\circle*{4}}

\put(280,70){$C_0$}
\put(85,130){$C_2$}
\put(175,100){$C_5$}
\put(215,120){$C_6$}

\put(84,70){$\tx_2$}
\put(174,70){$\tx_5$}
\put(214,70){$\tx_6$}

\put(55,70){$\tx_1$}
\put(115,70){$\tx_3$}
\put(140,70){$\tx_4$}
\put(250,70){$\tx_m$}

\put(230,70){$\cdots$}
\end{picture}
\caption{The curve $\BC$}
\end{figure}

Let $\tx_r=\br^{-1}(x_r)\cap C_0$ be the preimages of the points 
$x_r$ ($r=1,\dots,m$) on the component $C_0$ of $\BC$. Thus if 
$\theta_r\in \sE_{i-1}$ then $\tx_r$ is the point of intersection
of the components $C_0$ and $C_r$ of $\BC$.

\begin{lm} 
Let $\tx\in\BC$. Then we have
$$
\bq(\tx)\in\begin{cases}
K_{\lbr\theta_1,\dots,\theta_r\smile i,\dots,\theta_m\rbr}\subset K_{\alp}, &
\text{if $\tx\in C_r-\{\tx_r\}$},\\
K_{\lbr\theta_1,\dots,\theta_r,\dots,\theta_m,i\rbr}\subset K_{\alp}, &
\text{if $\tx\in C_0-\{\tx_1,\dots,\tx_m\}$},\\
K_{\lbr\theta_1,\dots,\lbr\theta_r,i\rbr,\dots,\theta_m\rbr}\subset K_{\alp}, &
\text{if $\tx=\tx_r$}.
\end{cases}
$$
\end{lm}
\begin{pf}
Easy.
\end{pf}

\begin{lm}
Consider the stratum $K_{\lbr\theta_1,\dots,\lbr\theta_r,i\rbr,\dots,\theta_m\rbr}$
of $K_{\alp}$. Then it lies in the component 
$\OK_{\lbr\theta_1,\dots,\theta_r,\dots,\theta_m,i\rbr}$ of $K_{\alp}$.
Other $(|\al|+1)$-dimensional components of $K_{\alp}$,
containing this stratum are listed below:
$$
\begin{array}{ll}
\!\!\OK_{\lbr\theta_1,\dots,\theta_r\smile i,\dots,\theta_m\rbr}\ \text{ and }\ 
\OK_{\lbr\theta_1,\dots,i\smile \theta_r,\dots,\theta_m\rbr} &
\!\text{if $\theta_r\in \sE_{i-1}\cap \sB_{i+1}$}\\
\!\!\OK_{\lbr\theta_1,\dots,\theta_r\smile i,\dots,\theta_m\rbr} &
\!\text{if $\theta_r\in \sE_{i-1}$}\\
\!\!\OK_{\lbr\theta_1,\dots,i\smile \theta_r,\dots,\theta_m\rbr} &
\!\text{if $\theta_r\in \sB_{i+1}$}\\
\end{array}
$$
This list is complete.
\end{lm}
\begin{pf}
Evident.
\end{pf}

Given a Kostant partition $A'\in\FK(\alp)$ consider the intersection
$$
\FE_A^{A'}=\FE_\al^i\cap\bp^{-1}(\OK_A)\cap\bq^{-1}(\OK_{A'}).
$$
We are interested in irreducible components of $\FE_A^{A'}$
which are dominant over $\OK_A$.

\begin{pr}\label{pdom1}
If $A'={\lbr\theta_1,\dots,\theta_r,\dots,\theta_m,i\rbr}$ then 
$\FE_A^{A'}$ has only one irreducible component, dominant
over $\OK_A$. This component is a generically $C$-fibration
over $\OK_A$. Its fiber over a generic point $E_\bullet\in K_A\subset\OK_A$
is equal to the component $C_0$ of the curve $\BC=\bp^{-1}(E_\bullet)$.
\end{pr}
\begin{pf}
Easy.
\end{pf}

\begin{pr}\label{pdom2}
If $\theta_r\in \sE_{i-1}$ and $A'=\lbr\theta_1,\dots,\theta_r\smile i,\dots,\theta_m\rbr$
then $\FE_A^{A'}$ has only one irreducible component, dominant
over $\OK_A$. This component is a generically 
$(\underbrace{\PP^1\sqcup\dots\sqcup\PP^1}_{\text{$m(\theta_r,A)$ times}})$-fibration
over $\OK_A$. Its fiber over a generic point $E_\bullet\in K_A\subset\OK_A$
is equal to the disjoint union of all components $C_{r'}$ of the curve 
$\BC=\bp^{-1}(E_\bullet)$ with $\theta_{r'}=\theta_r$.
\end{pr}
\begin{pf}
Easy.
\end{pf}

\begin{pr}\label{pdom3}
If $\theta_r\in \sB_{i+1}$ and $A'=\lbr\theta_1,\dots,i\smile \theta_r,\dots,\theta_m\rbr$
then $\FE_A^{A'}$ has only one irreducible component, dominant over $\OK_A$. 
This component is a generically finite $m(\theta_r,A)$-fold covering of $\OK_A$. 
Its fiber over a generic point $E_\bullet\in K_A\subset\OK_A$ is equal to the 
set all points $\tx_{r'}$ of the curve $\BC=\bp^{-1}(E_\bullet)$ 
such that $\theta_{r'}=\theta_r$.
\end{pr}
\begin{pf}
Easy.
\end{pf}

\section{$\hfsl_n$-module structure}

\subsection{Preliminaries}

We want to introduce a $\hfsl_n$-module structure on the vector space 
$$
\BM=\bigoplus_{\al\in\NN[I]}H^0(K_\al,\QQ).
$$
This space is naturally $Y$-graded. Recall that the affine Lie algebra $\hfsl_n$ 
is given by the generators $e_i,f_i,h_i$, $(i\in I)$ (Chevalley generators)
which satisfy the following set of relations (Serre relations)
$$
\begin{array}{ll}
\ad(e_i)^{(1-a_{ij})}e_j=0, & (i\ne j\in I)\\
\ad(f_i)^{(1-a_{ij})}f_j=0, & (i\ne j\in I)\\
{}[e_i,f_j]=\delta_{ij}h_i, & (i,j\in I)\\
{}[h_i,e_j]=a_{ij} e_j, & (i,j\in I)\\
{}[h_i,f_j]=-a_{ij} f_j, & (i,j\in I)\\
{}[h_i,h_j]=0, & (i,j\in I)
\end{array}
$$


Recall also, that $\bc=\sum_{i\in I}\hp_i$
is the central element of $\hfsl_n$.

\subsection{Definition of $e_i$}

We choose the natural $\QQ$-basis in the spaces $H^0(K_\al,\QQ)$,
formed by the fundamental classes of topdimensional irreducible
components of~$K_\al$. Let $v_A=[\OK_A]\in H^0(K_\al,\QQ)$ be
the fundamental class of the component $\OK_A$, corresponding to
a Kostant partition $A\in\FK(\al)$.

We will begin with the definition of  the action of the Chevalley 
generators $e_i$, $f_i$ and $\hp_i$ and after that we will check 
that the Serre relations between them are satisfied.

\begin{df}
The operator $e_i:H^0(K_\al,\QQ)\to H^0(K_{\alp},\QQ)$
is given by the correspondence $\FE_\al^i\subset K_\al\times K_{\alp}$. 
\end{df}

Let $\eps_i(A,A')$ denote the matrix coefficient of the
operator $e_i$ with respect to our basises. Thus for
all $A\in\FK(\al)$ we have
\begin{equation}\label{ei}
e_i(v_A)=\sum_{A'\in\FK(\alp)}\eps_i(A,A')v_{A'}.
\end{equation}

\begin{pr}\label{epsi}
Let $A'=\lbr\theta'_1,\dots,\theta'_m\rbr$.
We have
$$
\eps_i(A,A')=\begin{cases}
m(\theta_r,A'), &\text{if $\theta'_r\in \sE_i$ and 
$A=\lbr\theta'_1,\dots,\theta'_r\frown i,\dots,\theta'_m\rbr$}\\
0, &\text{otherwise}
\end{cases}
$$
\end{pr}
\begin{pf}
Follows from the Proposition \ref{qdom} and from the Lemma \ref{qdegree}.
\end{pf}

\subsection{Definition of $\hp_i$}

\begin{df}
We define the operator $\hp_i:H^0(K_\al,\QQ)\to H^0(K_\al,\QQ)$
as a scalar multiplication by $\langle i',(2-2g)\rho+\al_0+\al\rangle$.
\end{df}

\subsection{Definition of $f_i$}

The definition of the operator $f_i$ is rather more complicated.
The reason is that the dimension of $\FE_\al^i$ is equal to $|\al|+1$,
so it doesn't define an operator $H^0(K_{\alp},\QQ)\to H^0(K_\al,\QQ)$.
Thus, as in the case of excess intersection we should introduce certain
second cohomology classes $\xi_\al^i\in H^2(\FE_\al^i)$.

This can be done as follows.
Consider the space $\FE_\al^i\times S$ and let
$E_\bullet\subset\CO_{\FE_\al^i}\boxtimes\CE_\bullet$ 
(resp. $E'_\bullet\subset\CO_{\FE_\al^i}\boxtimes\CE_\bullet$) 
denote the universal degree-$\al$ (resp. degree-$(\alp)$)
periodic subflag on $\FE_\al^i\times S$. We have the universal
embedding $E'_\bullet\subset E_\bullet$ which gives  the 
exact sequence
$$
0 \to E'_i \to E_i \to (\id\times\bj)_*\Delta^i_*L_\al^i \to 0
$$
where $L_\al^i$ is a line bundle on $\FE_\al^i$ and
the embedding $\Delta^i$ is defined from the following 
cartesian square
$$
\begin{CD}
\FE_\al^i\times S @<{\id\times\bj}<< \FE_\al^i\times C @<{\Delta^i}<< \FE_\al^i\\
@V{\br\times\id}VV		    @V{\br\times\id}VV		     @V{\br}VV\\
C\times S	  @<{\id\times\bj}<< C\times C	      @<{\Delta}<<   C
\end{CD}
$$

Consider also the following commutative diagram
of sheaves on $K_\al\times S$
\begin{equation}\label{deffai}
\begin{CD}
  @.   0		@.   0	       @.   0	 \\
@.     @VVV		     @VVV	    @VVV \\
0 @>>> E_{i+1}\cap\CE_i @>>> E_{i+1}   @>>> E_{i+1}/(E_{i+1}\cap\CE_i) @>>> 0 \\
@.     @VVV		     @VVV	    @VVV \\
0 @>>> \CE_i		@>>> \CE_{i+1} @>>> (\id\times\bj)_*\CL_{i+1}   @>>> 0 \\
@.     @VVV		     @VVV	    @VVV \\
0 @>>> T'_i		@>>> T_{i+1}   @>>> T_{i+1}/T'_i	       @>>> 0 \\
@.     @VVV		     @VVV	    @VVV \\
  @.   0		@.   0	       @.   0	
\end{CD}
\end{equation}
where $T'_i=\CE_i/(E_{i+1}\cap\CE_i)=\Im(T_i\to T_{i+1})$.
It follows that $E_{i+1}/(E_{i+1}\cap\CE_i)$ is a subsheaf in 
$(\id\times\bj)_*\CL_{i+1}$, hence 
$$
E_{i+1}/(E_{i+1}\cap\CE_i)\cong(\id\times\bj)_*\CF_\al^i
$$
where $\CF_\al^i$ is a torsion free rank 1 sheaf on $K_\al\times C$.

\begin{df}
Let 
$$
\xi_\al^i=c_1((\Delta^i)^*(\bp\times\id)^*\CF_\al^i)-
c_1(L_\al^i)+c_1(\br^*\CT_C)\in H^2(\FE_\al^i,\ZZ),
$$
where $\CT_C$ stands for the tangent bundle of the curve $C$.
\end{df}
\begin{rm}
It is clear that first and second summands here depend not only on $i \bmod n$
but on $i$ itself also. However, it will be seen from the proof of the 
Proposition~\ref{phii} that the class $\xi_\al^i$ depend only on $i \bmod n$.
\end{rm}

\begin{df}
Let $A\in\FK(\al)$, $A'\in\FK(\alp)$ be a pair of Kostant partitions.
We define the operator $f_A^{A'}:H^0(\OK_{A'})\to H^0(\OK_A)$
as follows. If the intersection $\FE_A^{A'}$ is $(|\al|+1)$-dimensional,
then $f_A^{A'}$ is given by the class $(\xi_\al^i)|_{\FE_A^{A'}}$
(that is, $f_A^{A'}$ equals the integral of $\xi_\al^i$ over a generic
fiber of $\FE_A^{A'}$ over $K_A$).
If the intersection $\FE_A^{A'}$ is $|\al|$-dimensional,
then $f_A^{A'}$ is given by $\FE_A^{A'}$ (that is, 
$f_A^{A'}$ equals the cardinality of a generic
fiber of $\FE_A^{A'}$ over $K_A$).
 Finally, we define the
operator $f_i:H^0(K_{\alp},\QQ)\to H^0(K_\al,\QQ)$ as the
sum of operators $f_A^{A'}$ for all pairs $(A,A')$:
$$
f_i=\sum_{A\in\FK(\al),\ A'\in\FK(\alp)}f_A^{A'}.
$$
\end{df}

Let $\phi_i(A',A)$ denote the matrix coefficient of the
operator $f_i$ with respect to our bases. Thus for
all $A'\in\FK(\alp)$ we have
\begin{equation}\label{fi}
f_i(v_{A'})=\sum_{A\in\FK(\al)}\phi_i(A',A)v_A.
\end{equation}

\begin{pr}\label{phii}
Let $A=\lbr\theta_1,\dots,\theta_m\rbr$.
We have
$$
\phi_i(A',A)=\!\begin{cases}
M(i,A), &\!\!\!\!\text{if $A'=\lbr\theta_1,\dots,\theta_m,i\rbr$},\\
-m(\theta_r,A), &\!\!\!\!\text{if $\theta_r\in \sE_{i-1}$ and 
$A'=\lbr\theta_1,\dots,\theta_r\smile i,\dots,\theta_m\rbr$}\\
m(\theta_r,A), &\!\!\!\!\text{if $\theta_r\in \sB_{i+1}$ and 
$A'=\lbr\theta_1,\dots,i\smile \theta_r,\dots,\theta_m\rbr$}\\
0, &\!\!\!\!\text{otherwise}
\end{cases}
$$
where 
$$
M(i,A)=\langle i',(2-2g)\rho+\al_0\rangle+
\sum_{\theta\in \sB_{i+1}}(m(i\smile \theta,A)-m(\theta,A)).
$$
\end{pr}
\begin{pf}
The third case follows immediately from the Proposition~\ref{pdom3}
and from the definition of the operator $f_A^{A'}$. In order to check
the first two cases we should compute the restriction of the class $\xi_\al^i$
to the components $C_0$ and $C_r$ of the curve $\BC=\bp^{-1}(E_\bullet)$,
where $E_\bullet$ is a generic point of the stratum $K_A\subset K_\al$.

Recall that $\xi_\al^i$ is defined as a sum of three summands. We begin with
the computation of $((\Delta^i)^*(\bp\times\id)^*\CF_\al^i)_{|\BC}.$
Note that $(\bp\times\id)\cdot\Delta^i=\bp\times\br$, hence
$(\Delta^i)^*(\bp\times\id)^*\CF_\al^i=(\bp\times\br)^*\CF_\al^i$.
Since $\BC=\bp^{-1}(E_\bullet)$, we have the following cartesian square
$$
\begin{CD}
\BC				@>>>	\FE_\al^i		\\
@V{\br}VV				@V{\bp\times\br}VV	\\
\{E_\bullet\}\times C		@>>>	K_\al\times C
\end{CD}
$$
hence $((\bp\times\br)^*\CF_\al^i)_{|\BC}=
\br^*({\CF_\al^i}_{|\{E_\bullet\}\times C})$. 
It follows that this bundle is trivial on vertical components
of $\BC$ over $C$, that is
$$
c_1(((\Delta^i)^*(\bp\times\id)^*\CF_\al^i)_{|C_r})=0.
$$
At the same time its restriction to $C_0$ is isomorphic
to ${\CF_\al^i}_{|\{E_\bullet\}\times C}$, which according to the 
definition of $\CF_\al^i$ is equal to $E_{i+1}/(E_{i+1}\cap\CE_i)$.
The diagram \refeq{deffai} implies that
$E_{i+1}/(E_{i+1}\cap\CE_i)\cong \bj_*(\CL_{i+1}(-\dim\Ga(T_{i+1}/T'_i)))$.
So, finally we get
\begin{multline*}
c_1(((\Delta^i)^*(\bp\times\id)^*\CF_\al^i)_{|C_0})=
\deg\CL_{i+1}(-\dim\Ga(T_{i+1}/T'_i))=\\=
\deg\CL_{i+1}-\dim\Ga(T_{i+1}/T'_i)=
\deg\CL_{i+1}-\sum_{\theta\in \sB_{i+1}}m(\theta,A).
\end{multline*}

The second summand can be computed as follows. First consider 
a vertical component $C_r$. Note that $C_r$ is by definition
isomorphic to the projective line $\PP(H_r)$, where
$H_r=\Hom(E_\bullet,M_i\otimes\CO_{x_r})=\Hom(E_i/E_{i-1},\CO_{x_r}).$
Consider the natural morphism of 
coherent sheaves on $C_r\times S$
$$
\CO_{C_r}\boxtimes E_i \to 
\CO_{C_r}\boxtimes E_i/E_{i-1} \to 
\CO_{C_r}\boxtimes (H_r^*\otimes\CO_{x_r}) \to 
\CO_{C_r}(1)\boxtimes \CO_{x_r} 
%
%
$$
It is clear that the kernel of this morphism is nothing but
the restriction of the sheaf $E'_i$ from $\FE_\al^i\times S$
to $C_r\times S$. Hence ${L_\al^i}_{|C_r}\cong\CO_{C_r}(1)$
and 
$$
c_1({L_\al^i}_{|C_r})=1.
$$

Now we will treat the component $C_0$. Consider the following morphism 
of coherent sheaves on $C_0\times S$
\begin{multline*}
\CO_{C_0}\boxtimes E_i \to
\CO_{C_0}\boxtimes E_i/E_{i-1} \to
\CO_{C_0}\boxtimes E_i/(E_i\cap\CE_{i-1}) \to \\ \to
\CO_{C_0}\boxtimes \bj_*\CL_i(-\dim\Ga(T_i/T'_{i-1}))\to
(\id\times\bj)_*\CL_i(-\dim\Ga(T_i/T'_{i-1}))
%
%
\end{multline*}
It is clear that the kernel of this morphism is just
the restriction of the the sheaf $E'_i$ from $\FE_\al^i\times S$
to $C_0\times S$. Hence ${L_\al^i}_{|C_0}\cong\CL_i(-\dim\Ga(T_i/T'_{i-1}))$
and 
$$
c_1({L_\al^i}_{|C_0})=\deg\CL_i-\dim\Ga(T_i/T'_{i-1})=
\deg\CL_i-\sum_{\theta\in \sB_i}m(\theta,A).
$$

Finally, we should compute the third summand. But it is evident
that the bundle $\br^*\CT_C$ is trivial on the vertical components
of the curve $\BC$ and is isomorphic to $\CT_{C_0}$ on the 
component $C_0$. So summing up all contributions we obtain
\begin{multline*}
{\xi_\al^i}_{|C_0}=
\deg\CL_{i+1}-\sum_{\theta\in \sB_{i+1}}m(\theta,A)-
(\deg\CL_i-\sum_{\theta\in \sB_i}m(\theta,A))+2-2g=\\=
2-2g+\langle i',\al_0\rangle+
(\sum_{\theta\in \sB_i}m(\theta,A))-\sum_{\theta\in \sB_{i+1}}m(\theta,A))=M(i,A)
\end{multline*}
in the case of the component $C_0$ and
$$
{\xi_\al^i}_{|C_r}=0-1+0=-1
$$
in the case of the component $C_r$. So it remains to recall that
according to the Proposition~\ref{pdom2} the number of vertical
components, contributing to the coefficient $\phi_i(A',A)$
is equal to $m(\theta_r,A)$.
%
%
\end{pf}

We finish this section with the following.

\begin{th}\label{hfsln}
The operators $e_i,h_i,f_i$ $(i\in I)$ provide the vector space $\BM$
with a structure of $\hfsl_n$-module with the lowest weight $(2-2g)\rho+\al_0$
and the central charge $c_0=(2-2g)n+d$.
\end{th}

\section{Serre relations}

This section is devoted to the proof of the Theorem~\ref{hfsln}.

\subsection{Representation in differential operators}

Let $\BN=\QQ[x_\theta]$ denote a $Y$-graded vector space of polynomials
in variables $x_\theta$, $\theta\in R^+$, with the inverse grading
$$
\deg x^\theta=-\dim\theta
$$
Given a Kostant partition $A=\lbr\theta_1,\dots,\theta_m\rbr\in\FK(\al)$
let $x^A$ denote the following monomial
$$
x^A=x_{\theta_1}\cdot\dots\cdot x_{\theta_m}\in\BN_{-\al}.
$$
It is clear that monomials $x^A$ with $A\in\FK(\al)$ form a basis 
of $\BN_{-\al}$. We define the pairing $\BN\otimes\BM\to\QQ$ as follows
$$
\langle x^A,v_{A'}\rangle=\delta_{A,A'}.
$$
This gives an isomorphism $\BN\cong\BM^*$. Let 
$$
e_i^T:\BN_{-\al-i}\to\BN_{-\al},\quad
\hp_i^T:\BN_{-\al}\to\BN_{-\al},\quad
f_i^T:\BN_{-\al}\to\BN_{-\al-i},
$$
denote the adjoint operators of the operators $e_i$, $\hp_i$ and $f_i$,
defined in the previous section. 

\begin{lm}
We have 
$$
e_i^T(x^{A'})=\sum_{A\in\FK(\al)}\eps_i(A,A')x^A,\quad
f_i^T(x^A)=\sum_{A'\in\FK(\alp)}\phi_i(A',A)x^{A'},
$$
where $\eps_i(A,A')$ and $\phi_i(A',A)$ are given by
the Propositions~\ref{epsi} and \ref{phii}, and 
${\hp_i^T}_{|\BN_{-\al}}$ is the scalar multiplication
by $\langle i',(2-2g)\rho+\al_0+\al\rangle$.
\end{lm}
\begin{pf}
Evident.
\end{pf}

The following Proposition shows that the operators $e_i^T$, $f_i^T$ 
are in fact first order differential operators.

\begin{pr}\label{diff}
We have
\begin{equation}\label{ehf}
\rnc{\arraystretch}{1.5}
\begin{array}{lll}
\displaystyle
e_i^T & = & \sum\limits_{\theta\in \sE_{i-1}} x_\theta\dx_{\theta\smile i},\\
\hp_i^T & = & \sum\limits_{\theta\in \sE_i}x_\theta\dx_\theta - 
\sum\limits_{\theta\in \sE_{i-1}}x_\theta\dx_\theta +
\sum\limits_{\theta\in \sB_i}x_\theta\dx_\theta - 
\sum\limits_{\theta\in \sB_{i+1}}x_\theta\dx_\theta +c_i,\\
f_i^T & = & \sum\limits_{\theta\in \sB_{i+1}} x_{i\smile \theta}\dx_\theta - 
\sum\limits_{\theta\in \sE_{i-1}}x_{\theta\smile i}\dx_\theta + 
x_i(\sum\limits_{\theta\in \sB_i} x_{\theta}\dx_\theta - 
\sum\limits_{\theta\in \sB_{i+1}} x_{\theta}\dx_\theta +c_i),
\end{array}
\end{equation}
where $\dx_0=0$ and $c_i=\langle i',(2-2g)\rho+\al_0\rangle$.
\end{pr}
\begin{pf}
Direct calculations.
\end{pf}

It is clear that the Serre relations for the operators $e_i,h_i,f_i$
are equivalent to the Serre relations for the operators $e_i^T,h_i^T,f_i^T$.
Hence the Theorem~\ref{hfsln} follows from the following.

\begin{th}\label{key}
The operators \refeq{ehf} provide the vector space $\BN$ with
a structure of $\hfsl_n$-module with the lowest weight $c_i$.
\end{th}
The proof of the Theorem~\ref{key} will take the rest of the section.

\subsection{Check of the relations}

Now we can apply the Proposition~\ref{diff} to verify
the Serre relations. We begin with the following useful
notation.

Given a raiz $\theta\in \sB_i-\{0\}\subset R^+$ beginning at $i$
we define the differential operators
$$
\BE(\theta)=\sum_{\vartheta\in 
\sE_{i-1}}x_{\vartheta\smile\theta}\dx_\vartheta\quad\text{and}\quad
\TBE(\theta)=\sum_{\vartheta\in 
\sE_{i-1}}x_{\vartheta}\dx_{\vartheta\smile\theta}.
$$
Similarly, given a raiz $\theta\in \sE_i-\{0\}\subset R^+$ ending at $i$ we define
$$
\BB(\theta)=\sum_{\vartheta\in \sB_{i+1}}x_{\theta\smile\vartheta}\dx_\vartheta.
$$
Finally we define 
$$
\BE_i=\sum_{\theta\in \sE_{i-1}}x_{\theta}\dx_\theta,\qquad
\BB_i=\sum_{\theta\in \sB_{i+1}}x_{\theta}\dx_\theta
$$
and 
$$
\Delta_i=\BB_{i-1}+\BB_i+c_i.
$$
Using this notation we can write our operators in a more compact form
\begin{equation}\label{compactform}
\rnc{\arraystretch}{1.2}
\begin{array}{lll}
e_i^T & = & \TBE(i),\\
\hp_i^T & = & \BE_{i+1} - \BE_i + \Delta_i ,\\
f_i^T & = & \BB(i)-\BE(i)+x_i\Delta_i.
\end{array}
\end{equation}

We will make also the following agreement: 
$$
\text{if $\theta\not\in R_+$ then $\BE(\theta)=\TBE(\theta)=\BB(\theta)=x_\theta=0$.}
$$

\begin{lm}\label{comm1}
We have
$$
\rnc{\arraystretch}{1.2}
\begin{array}{lll}
{}[\TBE(\theta_1),\TBE(\theta_2)] & = & \TBE(\theta_1\smile\theta_2)-\TBE(\theta_2\smile\theta_1),\\
{}[\BE(\theta_1),\BE(\theta_2)] & = & \BE(\theta_2\smile\theta_1)-\BE(\theta_1\smile\theta_2),\\
{}[\BB(\theta_1),\BB(\theta_2)] & = & \BB(\theta_1\smile\theta_2)-\BB(\theta_2\smile\theta_1),\\
{}[\BB_i,\BB_j] & = & 0,\\
{}[\BB(\theta),\BB_i] & = & \begin{cases}
\BB(\theta), &\text{if $\theta\in \sE_i\setminus \sB_{i+1}$,}\\
-\BB(\theta), &\text{if $\theta\in \sB_{i+1}\setminus \sE_i$,}\\
0, &\text{otherwise.}
\end{cases}
\end{array}
$$
\end{lm}
\begin{pf}
Easy.
\end{pf}

\begin{lm}\label{comm2}
We have
$$
\rnc{\arraystretch}{1.2}
\begin{array}{lll}
{}[\TBE(i),\BB(j)] & = & 0,\\
{}[\TBE(i),\BB_j] & = & 0,\\
{}[\TBE(i),\BE(j)] & = & \delta_{ij}(\BE_i-\BE_{i+1}),\\
{}[\TBE(i),x_j] & = & \delta_{ij}.
\end{array}
$$
\end{lm}
\begin{pf}
Easy.
\end{pf}

\begin{lm}\label{comm3}
We have
$$
\rnc{\arraystretch}{1.2}
\begin{array}{lll}
{}[\BE(\theta_1),\BB(\theta_2)] & = & 0,\\
{}[\BE(\theta),\BB_i] & = & 0,\\
{}[\BE(\theta_1),x_{\theta_2}] & = & x_{\theta_2\smile\theta_1},\\
{}[\BB(\theta_1),x_{\theta_2}] & = & x_{\theta_1\smile\theta_2},\\
{}[\BB_i,x_{\theta}] & = & \begin{cases}
x_{\theta}, &\text{if $\theta\in \sB_{i+1}$,}\\
0, &\text{otherwise.}
\end{cases}
\end{array}
$$
\end{lm}
\begin{pf}
Easy.
\end{pf}


\begin{lm}\label{comm4}
We have 
$$
\rnc{\arraystretch}{1.2}
\begin{array}{lll}
{}[\Delta_i,\BE(\theta)] & = & 0;\\
{}[\Delta_i,\TBE(\theta)] & = & 0;\\
{}[\Delta_i,\BB(\theta)] & = & \langle i',\dim\theta\rangle\BB(\theta);\\
{}[\Delta_i,x_\theta] & = & \begin{cases}
x_{\theta}, &\text{if $\theta\in \sB_i$,}\\
-x_{\theta}, &\text{if $\theta\in \sB_{i+1}$,}\\
0, &\text{otherwise};
\end{cases}\\
{}[\Delta_i,\Delta_j] & = & 0.
\end{array}
$$
\end{lm}
\begin{pf}
Easy.
\end{pf}

Now we are ready to check the relations.

\subsection{Serre relations for $e_i^T$}

\begin{pr}
If $j\ne i\pm1$ then $[e_i^T,e_j^T]=0$.

\noindent
If $j=i\pm1$ and $n\ne 2$ then $\ad(e_i^T)^2e_j^T=0$.

\noindent
If $j=i+1$ and $n=2$ then $\ad(e_i^T)^3e_j^T=0$.
\end{pr}
\begin{pf}
In the first case $i\smile j\not\in R^+$ and $j\smile i\not\in R^+$ hence according to \refeq{compactform}
and to the Lemma~\ref{comm1} we get 
$$
[\TBE(i),\TBE(j)]=0.
$$
If $j=i+1$ and $n\ne 2$ then 
$$
[\TBE(i),[\TBE(i),\TBE(j)]]=[\TBE(i),\TBE(i\smile j)]=0
$$
and if $j=i-1$ and $n\ne 2$ then 
$$
[\TBE(i),[\TBE(i),\TBE(j)]]=-[\TBE(i),\TBE(j\smile i)]=0.
$$
Finally, if $j=i+1$ and $n=2$ then
\begin{multline*}
[\TBE(i),[\TBE(i),[\TBE(i),\TBE(j)]]]=[\TBE(i),[\TBE(i),\TBE(i\smile j)-\TBE(j\smile i)]]=\\=
-2[\TBE(i),\TBE(i\smile j\smile i)]=0.
\end{multline*}
\end{pf}

\subsection{Commutators of $e_i^T$ and $f_j^T$}

\begin{pr}
We have $[e_i^T,f_j^T]=\delta_{ij}\hp_i^T$.
\end{pr}
\begin{pf}
According to \refeq{compactform}
and to the Lemma \ref{comm2} we have
$$
[\TBE(i),\BB(j)-\BE(j)+x_j\Delta_j]=\\=
-\delta_{ij}(\BE_i-\BE_{i+1})+\delta_{ij}\Delta_j=\delta_{ij}\hp_i^T.
$$
\end{pf}

\subsection{Serre relations for $f_i^T$}

\begin{pr}
If $j\ne i\pm1$ then $[f_i^T,f_j^T]=0$.
\end{pr}
\begin{pf}
According to \refeq{compactform}
and to the Lemmas~\ref{comm1}, \ref{comm3} and \ref{comm4} we have
$$
[f_i^T,f_j^T]=[\BB(i)-\BE(i)+x_i\Delta_i,\BB(j)-\BE(j)+x_j\Delta_j]=0.
$$
\end{pf}

\begin{pr}
If $j=i\pm1$ and $n\ne2$ then $\ad(f_i^T)^2f_j^T=0$.
\end{pr}
\begin{pf}
Assume that $j=i+1$. We denote 
$$
D=[f_i^T,f_j^T]=[\BB(i)-\BE(i)+x_i\Delta_i,\BB(j)-\BE(j)+x_j\Delta_j].
$$
Then according to \refeq{compactform}
and to the Lemmas~\ref{comm1}, \ref{comm3} and \ref{comm4} we have
$$
D=\BB(i\smile j)-\BE(i\smile j)+x_j\BB(i)-x_i\BB(j)+x_{i\smile j}\Delta_i+(x_{i\smile j}-x_ix_j)\Delta_j.
$$
It suffices to show that $[f_i^T,D]=[f_j^T,D]=0$.
Aplying once more \refeq{compactform} and the Lemmas~\ref{comm1}, \ref{comm3} and 
\ref{comm4} we get
$$
\rnc{\arraystretch}{1.2}
\arraycolsep=2pt
\begin{array}{lll}
{}[\BB(i),D] & = & 
x_{i\smile j}\BB(i)-x_i\BB(i\smile j)-2x_{i\smile j}\BB(i)+(x_{i\smile j}-x_ix_j)\BB(i)-x_ix_{i\smile j}\Delta_j;\\
{}[\BE(i),D] & = & 0;\\
{}[x_i\Delta_i,D] & = &
x_i\BB(i\smile j)+x_ix_j\BB(i)+x_ix_{i\smile j}\Delta_j.
\end{array}
$$
Hence $[f_i^T,D]=0$. Similarly, we get
$$
\rnc{\arraystretch}{1.2}
\begin{array}{lll}
{}[\BB(j),D] & = & -x_j\BB(i\smile j)+x_{i\smile j}\BB(j)-2(x_{i\smile j}-x_ix_j)\BB(j);\\
{}[\BE(j),D] & = & -x_{i\smile j}\BB(j)-x_{i\smile j}x_j\Delta_j;\\
{}[x_j\Delta_j,D] & = & 
x_j\BB(i\smile j)-x_jx_{i\smile j}\Delta_j-2x_ix_j\BB(j)+x_jx_{i\smile j}\Delta_j-x_jx_{i\smile j}\Delta_j.
\end{array}
$$
Hence $[f_j^T,D]=0$. 
\end{pf}

\begin{pr}
If $j=i+1$ and $n=2$ then $\ad(f_i^T)^3f_j^T=0$.
\end{pr}
\begin{pf}
We denote 
$$
D=[f_i^T,f_j^T]=[\BB(i)-\BE(i)+x_i\Delta_i,\BB(j)-\BE(j)+x_j\Delta_j].
$$
Then according to \refeq{compactform}
and to the Lemmas~\ref{comm1}, \ref{comm3} and \ref{comm4} we have
\begin{multline*}
D=\BB(i\smile j)-\BB(j\smile i)-\BE(i\smile j)+\BE(j\smile i)+\\
+2x_j\BB(i)-2x_i\BB(j)+(x_{i\smile j}-x_{j\smile i}+x_ix_j)\Delta_i+(x_{i\smile j}-x_{j\smile i}-x_ix_j)\Delta_j.
\end{multline*}
Aplying once more \refeq{compactform} and the Lemmas~\ref{comm1}, \ref{comm3} 
and \ref{comm4} we get
$$
\rnc{\arraystretch}{1.2}
\begin{array}{lll}
{}[\BB(i),D] & = & 
-2\BB(i\smile j\smile i)+2x_{i\smile j}\BB(i)-2x_i\BB(i\smile j)+2x_i\BB(j\smile i)+\\
&& +(x_ix_{i\smile j}-x_{i\smile j\smile i})\Delta_i-2(x_{i\smile j}-x_{j\smile i}+x_ix_j)\BB(i)-\\
&& -(x_ix_{i\smile j}+x_{i\smile j\smile i})\Delta_j+2(x_{i\smile j}-x_{j\smile i}-x_ix_j)\BB(i);\\
{}[\BE(i),D] & = & -2\BE(i\smile j\smile i)+2x_{j\smile i}\BB(i)+\\
&& +(x_{i\smile j\smile i}+x_ix_{j\smile i})\Delta_i+(x_{i\smile j\smile i}-x_ix_{j\smile i})\Delta_j;\\
{}[x_i\Delta_i,D] & = & -2x_{i\smile j\smile i}\Delta_i+2x_ix_j\BB(i)+2x_i^2\BB(j)+2x_ix_{j\smile i}\Delta_i+\\
&& +x_i(2x_{i\smile j}+x_ix_j)\Delta_i+x_i(x_{i\smile j}+x_{j\smile i})\Delta_j-\\
&& -x_i(x_{i\smile j}-x_{j\smile i}-x_ix_j)\Delta_i.
\end{array}
$$
Hence $[f_i^T,D]=2\TD$, where
\begin{multline*}
\TD=-\BB(i\smile j\smile i)+\BE(i\smile j\smile i)+x_i\BB(j\smile i)-x_i\BB(i\smile j)+(x_{i\smile j}-x_{j\smile i}-x_ix_j)\BB(i)+\\
+x_i^2\BB(j)+(-2x_{i\smile j\smile i}+x_ix_{i\smile j}+x_ix_{j\smile i}-x_i^2x_j)\Delta_i+(-x_{i\smile j\smile i}+x_ix_{j\smile i})\Delta_j.
\end{multline*} 
Aplying once more \refeq{compactform} and the Lemmas~\ref{comm1}, \ref{comm3} 
and \ref{comm4} we get
$$
\arraycolsep=2pt
\rnc{\arraystretch}{1.2}
\begin{array}{lll}
{}[\BB(i),\TD] & = & 
2x_i\BB(i\smile j\smile i)-(x_{i\smile j\smile i}+x_ix_{i\smile j})\BB(i)+\\
&& +x_i^2(\BB(i\smile j)-\BB(j\smile i))+(x_ix_{i\smile j\smile i}-x_i^2x_{i\smile j})\Delta_i+\\
&& +2(2x_{i\smile j\smile i}-x_ix_{i\smile j}-x_ix_{j\smile i}+x_i^2x_j)\BB(i)+\\
&& +x_ix_{i\smile j\smile i}\Delta_j+2(-x_{i\smile j\smile i}+x_ix_{j\smile i})\BB(i);\\
{}[\BE(i),\TD] & = & (x_{i\smile j\smile i}-x_ix_{j\smile i})\BB(i)+(x_ix_{i\smile j\smile i}-x_i^2x_{j\smile i})\Delta_i;\\
{}[x_i\Delta_i,\TD] & = & -2x_i\BB(i\smile j\smile i)+x_i^2\BB(j\smile i)-x_i^2\BB(i\smile j)+x_ix_{i\smile j\smile i}\Delta_i+\\
&& +x_i(x_{i\smile j}+x_{j\smile i})\BB(i)-2x_i(-x_{i\smile j}+x_{j\smile i}+x_ix_j)\BB(i)-\\
&& -x_i^2x_{j\smile i}\Delta_i+x_i(-2x_{i\smile j\smile i}+2x_ix_{i\smile j}-x_i^2x_j)\Delta_i+\\
&& +x_i(2x_{i\smile j\smile i}-x_ix_{i\smile j}-x_ix_{j\smile i}+x_i^2x_j)\Delta_i-\\
&& -x_ix_{i\smile j\smile i}\Delta_j+x_i(-x_{i\smile j\smile i}+x_ix_{j\smile i})\Delta_i.
\end{array}
$$
Hence $[f_i^T,\TD]=0$ and the Proposition follows.
\end{pf}

\section{Extending the action of $\hfsl_n$ on $\BM$ to $\hfgl_n$}

In this section we will deal with the Lie algebras $\hfsl_n$ 
for various $n$. So in order to avoid a confusion
we will denote the corresponding set of simple raiz by $I(n)$, 
the system of raiz by $R^+(n)$, and the raiz lattice by $Y(n)$.

\subsection{$\hfsl_n$ and $\hfsl_{kn}$}

Recall that the group $\ZZ/n\ZZ$ acts on the
Lie algebra $\hfsl_n$ by outer automorphisms. This group acts also on the
set of simple raiz $I(n)$, on the raiz system $R^+(n)$, and 
on the raiz lattice $Y(n)$. We denote the action of an
element $a\in\ZZ/n\ZZ$ by $\tau_a$.

\begin{lm}\label{kinv}
For all integers $k\ge2$, $n\ge2$ there is a Lie algebras homomorphism 
$\mu:\hfsl_n\to\hfsl_{kn}$ defined on the Chevalley generators
as follows
$$
\mu(e_i)=\sum_{a\in n\ZZ/(kn)\ZZ}e_{\tau_a(i)},\quad
\mu(f_i)=\sum_{a\in n\ZZ/(kn)\ZZ}f_{\tau_a(i)},\quad
\mu(\hp_i)=\sum_{a\in n\ZZ/(kn)\ZZ}\hp_{\tau_a(i)}.
$$
\end{lm}
\begin{pf}
Evident.
\end{pf}

On the other hand, the identifications
$$
R^+(kn)=(\ZZ\letimes\ZZ)/(kn)\ZZ,\qquad
R^+(n)=(\ZZ\letimes\ZZ)/n\ZZ
$$
give rise to the projection
$$
\zeta:R^+(kn)\to R^+(n),\qquad(p,q) \bmod kn\mapsto(p,q)\bmod n.
$$
The projection $\zeta$ in its turn induces a morphism of algebras
$$
\zeta:\BN(kn)\to\BN(n),\qquad x_\vartheta\mapsto x_{\zeta(\vartheta)}.
$$

The following Proposition will be very important below.

\begin{pr}
The homomorphism $\zeta$ is a homomorphism of $\hfsl_n$-modules,
that is for all $\xi\in\hfsl_n$, $P(x)\in\BN(kn)$ we have
\begin{equation}\label{xi}
\xi\cdot\zeta(P)=\zeta(\mu(\xi)\cdot P).
\end{equation}
where $\cdot$ stands for the action of the Lie algebra $\hfsl_n$ 
$($resp. $\hfsl_{kn})$ on the space $\BN(kn)$ $($resp. $\BN(n))$.
\end{pr}
\begin{pf}
Since $\mu$ is a Lie algebra homomorphism it suffices to check~\refeq{xi}
only for the Chevalley generators $e_i$ and $f_i$. On the other hand,
both $e_i$ and $f_i$ act as first order differential operators, hence
both LHS and RHS of~\refeq{xi} are first order differential operators
on $\BN(kn)$ with values in $\BN(n)$. Hence it suffices to check~\refeq{xi}
only for $P=x_\vartheta$, $\vartheta\in R^+(kn)\cup\{0\}$, and this can be done
straightforwardly.
\end{pf}

In fact, the homomorphism $\mu:\hfsl_n\to\hfsl_{kn}\subset\Diff(\BN(kn))$
can be extended to some bigger subalgebra of $\Diff(\BN(n))$. 

\begin{df}
We define
$$
\mu(\TBE(\theta))=\sum_{\vartheta\in\zeta^{-1}(\theta)}\TBE(\vartheta),
\qquad\theta\in R^+(n).
$$
\end{df}

The analog of the property~\refeq{xi} is satisfied for the 
operators $\xi=\TBE(\theta)$:

\begin{pr}
We have
\begin{equation}\label{xi1}
\TBE(\theta)\cdot\zeta(P)=\zeta(\mu(\TBE(\theta))\cdot P).
\end{equation}
\end{pr}
\begin{pf}
Again it suffices to check~\refeq{xi1} only for $P=x_\vartheta$, 
$\vartheta\in R^+(kn)\cup\{0\}$ and this can be done straightforwardly.
\end{pf}

\subsection{The polynomial $P_n$}


Let $\alpha_n=\sum_{i\in I(n)}i\in Y(n)$
be the element of the lattice $Y(n)$, corresponding 
to the sequence $(\dots,1,1,1,1,\dots)$ of integers.

Consider the following element of the space $\BN(n)$
$$
P_n=\sum_{\ka=\lbr\theta_1\dots,\theta_m\rbr\in\FK(\alpha_n)}P_n^\ka x^\ka=
\sum_{\ka=\lbr\theta_1\dots,\theta_m\rbr\in\FK(\alpha_n)}
(-1)^{\sK(\ka)+1}x_{\theta_1}\cdot\ldots\cdot x_{\theta_m}.
$$

This element of the space $\BN(n)$ has the following properties.

\begin{lm}
Let $\theta\in R^+(n)$. We have
$$
\TBE(\theta)\cdot P_n=\begin{cases}
1, & \text{if $\dim\theta=\alpha_n$}\\
0, & \text{otherwise}
\end{cases}
$$
\end{lm}
\begin{pf}
The first case is evident, so assume that $\dim\theta\ne\alpha_n$.
It is clear that 
$$
\TBE(\theta)\cdot P_n=\sum_{\ka\in\FK(\alpha_n-\dim\theta)}\lambda_\ka x^\ka
$$
with some coefficients $\lambda_\ka\in\ZZ$. 
Hence we can assume also that $\alpha_n-\dim\theta\in\NN[I]$.

It suffices to show that all coefficients $\lambda_\ka$ vanish. 
Fix some $\ka=\lbr\theta_1,\dots,\theta_m\rbr$ and assume that $\theta\in \sB_i$. 
It is clear that the partition $\ka$ has exactly one element 
contained in $\sE_{i-1}$. Assume that it is $\theta_m$. 

Let $\ka'=\lbr\theta_1,\dots,\theta_m\smile\theta\rbr$ and 
$\ka''=\lbr\theta_1,\dots,\theta_m,\theta\rbr$.
It is clear that  
$$
\lambda_\ka=P_n^{\ka'}+P_n^{\ka''},
$$
but $\sK(\ka'')=\sK(\ka')+1$, hence $\lambda_\ka=0$.
\end{pf}

Let $f'_i=f_i^T-c_ix_i=\BB(i)-\BE(i)+x_i(\BB_{i-1}-\BB_i)$.

\begin{lm}
We have
$$
f'_i\cdot P_n=0
$$
for all $i\in I(n)$.
\end{lm}
\begin{pf}
It is clear that
$$
f'_i\cdot P_n=\sum_{\ka\in\FK(\alpha_n+i)}\lambda_\ka x^\ka
$$
with some coefficients $\lambda_\ka\in\ZZ$. We will check
that all coefficients $\lambda_\ka$ vanish.
We have the following cases to consider:

1) $\ka=\lbr i,i\rbr\cup\ka'$;

2) $\ka=\lbr i,i\smile\theta\rbr\cup\ka'$, $\theta\in \sB_{i+1}$;

3) $\ka=\lbr i,\theta\smile i\rbr\cup\ka'$, $\theta\in \sE_{i-1}$; 

4) $\ka=\lbr i,\theta\rbr\cup\ka'$, $\theta\ge i$ but $\theta\not\in \sB_i\cup \sE_i$;

5) $\ka=\lbr i\smile \theta_1,\theta_2\smile i\rbr\cup\ka'$, $\theta_1\in \sB_{i+1}$, 
$\theta_2\in \sE_{i-1}$.

In the first case we have
$$
\lambda_\ka=
P_n^{\lbr i\rbr\cup\ka'}((\BB_{i-1}-\BB_i)(x^{\lbr i\rbr\cup\ka'}))=0.
$$

In the second case we have
\begin{multline*}
\lambda_\ka=
P_n^{\lbr i,\theta\rbr\cup\ka'}+
P_n^{\lbr i\smile \theta\rbr\cup\ka'}((\BB_{i-1}-\BB_i)(x^{\lbr i\smile \theta\rbr\cup\ka'}))=\\=
P_n^{\lbr i,\theta\rbr\cup\ka'}+P_n^{\lbr i\smile\theta\rbr\cup\ka'}
\end{multline*}
which is zero because $\sK(\lbr i,\theta\rbr\cup\ka')=\sK(\lbr i\smile \theta\rbr\cup\ka')+1$.

In the third case we have
\begin{multline*}
\lambda_\ka=
-P_n^{\lbr i,\theta\rbr\cup\ka'}+
P_n^{\lbr \theta\smile i\rbr\cup\ka'}((\BB_{i-1}-\BB_i)(x^{\lbr \theta\smile i\rbr\cup\ka'}))=\\=
-P_n^{\lbr i,\theta\rbr\cup\ka'}-P_n^{\lbr i\smile \theta\rbr\cup\ka'}
\end{multline*}
which is zero because $\sK(\lbr i,\theta\rbr\cup\ka')=\sK(\lbr i\smile \theta\rbr\cup\ka')+1$.

In the fourth case we have
$$
\lambda_\ka=
P_n^{\lbr \theta\rbr\cup\ka'}((\BB_{i-1}-\BB_i)(x^{\lbr \ka\rbr\cup\ka'}))=0.
$$

Finally, in the fifth case we have
$$
\lambda_\ka=
P_n^{\lbr \theta_1,\theta_2\smile i\rbr\cup\ka'}-
P_n^{\lbr i\smile \theta_1,\theta_2\rbr\cup\ka'},
$$
which is zero because 
$\sK(\lbr \theta_1,\theta_2\smile i\rbr\cup\ka')=\sK(\lbr i\smile \theta_1,\theta_2\rbr\cup\ka')$.
\end{pf}


\subsection{$\hfgl_n$-structure}

Recall that 
$$
\hfgl_n/\bc=\hfsl_n/\bc\oplus\Heis/\bc,
$$ 
where $\bc$ is the central element and
$\Heis$ is the Heisenberg algebra, that is
$$
\Heis=\QQ(\dots,\ba_{-2},\ba_{-1},\ba_0,\ba_1,\ba_2,\dots)\oplus\QQ\bc
$$
with
$$
[\ba_p,\ba_q]=\delta_{p,-q}p\bc,\qquad [\ba_p,\bc]=0.
$$ 
Thus, if we want to extend an action of $\hfsl_n$ with the central charge
$c_0$ to the action of $\hfgl_n$, we have to construct an action of the 
Heisenberg algebra with the central charge~$nc_0$, commuting with the 
action of $\hfsl_n$.

\begin{df}
(compare ~\cite{s}) We define
\begin{equation}\label{heis}
\ba_p=\sum_{\dim\theta=p\alpha_n}\TBE(\theta),\quad
\ba_{-p}=c_0\zeta(P_{pn}),\quad p>0,\quad\text{and}\quad
\ba_0=c_0\id,\quad\bc=nc_0\id.
\end{equation}
\end{df}

\begin{pr}
The operators $\ba_p$ and $\bc$ satisfy the relations 
of the Heisenberg algebra.
\end{pr}
\begin{pf}
The equality $[\ba_p,\ba_q]=0$ for $p,q<0$ is evident and for
$p,q>0$ it follows from the Lemma~\ref{comm1}. Hence it suffices
to consider the case $p>0$, $q<0$. In this case we have
\begin{multline*}
[\ba_p,\ba_q]=
c_0\sum_{\dim\theta=p\alpha_n}\TBE(\theta)\cdot\zeta(P_{-qn})=
c_0\sum_{\dim\theta=p\alpha_n}\zeta(\mu(\TBE(\theta))\cdot P_{-qn})=\\=
c_0\sum_{\dim\zeta(\vartheta)=p\alpha_n}
\zeta(\TBE(\vartheta)\cdot P_{-qn})=
\begin{cases}
pnc_0, & \text{if $p=-q$}\\
0, & \text{otherwise}
\end{cases}
\end{multline*}
because $\dim\zeta(\vartheta)=p\alpha_n$ is equivalent to $\dim\vartheta=\alpha_{pn}$.
\end{pf}

Thus it remains to check that $\ba_p$ commute with the action of $\hfsl_n$.
We will need the following Lemma.

\begin{lm}\label{commagain}
We have
$$
\rnc{\arraystretch}{1.2}
\begin{array}{lll}
{}\sum_{\dim\theta=p\alpha_n}[\TBE(\theta),\BB(i)] & = &
\dx_{\theta_1}+x_i\dx_{\theta_2},\\
{}\sum_{\dim\theta=p\alpha_n}[\TBE(\theta),\BE(i)] & = &
\dx_{\theta_1}+x_i\dx_{\theta_3},\\
{}\sum_{\dim\theta=p\alpha_n}[\TBE(\theta),x_i\Delta_i] & = & 
x_i(\dx_{\theta_3}-\dx_{\theta_2}),
\end{array}
$$
where $\theta_1,\theta_2,\theta_3\in R^+(n)$ are defined by the following
properties
$$
\dim\theta_1=p\alpha_n-i,\qquad
\dim\theta_2=p\alpha_n\text{ and }\theta_2\in\sB_{i+1},\qquad
\dim\theta_3=p\alpha_n\text{ and }\theta_3\in\sB_i.
%
$$
\end{lm}
\begin{pf}
Evident.
\end{pf}

\begin{pr}
The operator $\ba_p$ commutes with $\hfsl_n$ for all $p>0$.
\end{pr}
\begin{pf}
It suffices to check that the operators $\ba_p$ commute with 
Chevalley generators $e_i^T$ and $f_i^T$. It follows from the Lemma~\ref{comm1} that
$$
[\ba_p,e_i^T]=\sum_{\dim\theta=p\alpha_n}[\TBE(\theta),\TBE(i)]=0 
$$
and the Lemma \ref{commagain} implies that $[\ba_p,f_i^T]=0$.
\end{pf}

\begin{pr}
The operator $\ba_{-p}$ commutes with $\hfsl_n$ for all $p>0$.
\end{pr}
\begin{pf}
We have
\begin{multline*}
[e_i^T,\ba_{-p}]=c_0 e_i^T\cdot\zeta(P_{pn})=
c_0\zeta(\mu(e_i^T)\cdot P_{pn})=\\=
c_0\sum\zeta(e_{\tau_a(i)}^T\cdot P_{pn})=
c_0\sum\zeta(\TBE(\tau_a(i))\cdot P_{pn})=0.
\end{multline*}
On the other hand
$$
[f_i^T,\ba_{-p}]=c_0 f_i\cdot\zeta(P_{pn})-c_0 c_ix_i\zeta(P_{pn})
$$
and
\begin{multline*}
f_i^T\cdot\zeta(P_{pn})=
\zeta(\mu(f^T_i)\cdot P_{pn})=\sum\zeta(f^T_{\tau_a(i)}\cdot P_{pn})=\\=
\sum\zeta((f'_{\tau_a(i)}+c_ix_{\tau_a(i)})\cdot P_{pn})=
\sum c_i\zeta(x_{\tau_a(i)})\zeta(P_{pn})=
c_ix_i\zeta(P_{pn})
\end{multline*}
and the Proposition follows.
\end{pf}

Thus we have proved the following.

\begin{th}\label{hfgln}
The operators \refeq{heis} extend the structure of $\hfsl_n$-module
of the vector space $\BN$ to that of $\hfgl_n$-module.
\end{th}

\begin{cor}
The vector space $\BM$ has the natural structure of a $\hfgl_n$-module.
\end{cor}

\begin{rem}
\label{meaning}
It is easy to show that the action of the operators $\ba_p$ with $p>0$
on the vector space $\BM$ can be described geometrically. Namely, the operator
$\ba_p$ is given by the correspondence
$$
\FE_\al^{p\alpha_n}=\{(E_\bullet,E'_\bullet)\in K_\al\times K_{\al+p\alpha_n}\ |\ 
\text{$E'_\bullet\subset E_\bullet$ and $\supp(E'_\bullet/E_\bullet)=\{x\}\in C$}\}
$$
in the same way as the operator $e_i$ is given by the correspondence $\FE_\al^i$.
It is tempting also to conjecture that the operator $\ba_{-p}$
is given by the same correspondence $\FE_\al^{p\alpha_n}$ in the same
way as the operator $f_i$ is given by the correspondence $\FE_\al^i$.
However, it is rather difficult to check, because the correspondence 
$\FE_\al^{p\alpha_n}$ has a lot of components with dimensions in the 
range $|\al|,\dots,|\al+p\alpha_n|$ and all of them should contribute 
to this operator.
\end{rem}

\end{document}